\documentclass[12pt]{amsart}
\usepackage{amsmath,amscd,amssymb,amsfonts,graphics}
\newcommand{\h}{\hbox}
\newcommand{\q}{\quad}
\newcommand{\nin}{\noindent}
\newcommand{\bs}{\par\bigskip}
\newcommand{\ms}{\par\medskip}
\newcommand{\sk}{\par\smallskip}
\newcommand{\bsn}{\par\bigskip\noindent}
\newcommand{\msn}{\par\medskip\noindent}
\newcommand{\skn}{\par\smallskip\noindent}
\newcommand{\ssb}{\raise.15ex\h{${\scriptscriptstyle\bullet}$}}
\newcommand{\ssc}{\,\raise.15ex\h{${\scriptstyle\circ}$}\,}
\newcommand{\msum}{\hbox{$\sum$}}
\newcommand{\mcup}{\hbox{$\bigcup$}}
\newcommand{\mopl}{\hbox{$\bigoplus$}}
\newcommand{\C}{{\mathbb C}}
\newcommand{\D}{{\mathbb D}}
\newcommand{\N}{{\mathbb N}}
\newcommand{\PP}{{\mathbb P}}
\newcommand{\Q}{{\mathbb Q}}
\newcommand{\RR}{{\mathbf R}}
\newcommand{\Z}{{\mathbb Z}}

\newcommand{\Hc}{{\mathcal H}}

\newcommand{\LL}{{\mathcal L}}

\newcommand{\OO}{{\mathcal O}}
\newcommand{\R}{{\mathcal R}}
\newcommand{\alt}{\widetilde{\al}}
\newcommand{\at}{\widetilde{a}}
\newcommand{\bt}{\widetilde{b}}
\newcommand{\ct}{\widetilde{c}}
\newcommand{\ddt}{\widetilde{d}}
\newcommand{\Rt}{\widetilde{\R}}
\newcommand{\Ht}{\widetilde{H}}

\newcommand{\mut}{\widetilde{\mu}}
\newcommand{\nut}{\widetilde{\nu}}
\newcommand{\rhot}{\widetilde{\rho}}
\newcommand{\chit}{\widetilde{\chi}}
\newcommand{\chih}{\widehat{\chi}}
\newcommand{\bo}{\overline{b}}
\newcommand{\cho}{\overline{\chi}}
\newcommand{\ddd}{{\rm d}}
\newcommand{\dd}{\partial}
\newcommand{\df}{{\rm d}f}

\newcommand{\ee}{{\bf e}}
\newcommand{\mm}{{\mathfrak m}}

\newcommand{\la}{\lambda}
\newcommand{\al}{\alpha}
\newcommand{\Ff}{F_{\!f}}

\newcommand{\kods}{\frac{k}{d}}
\newcommand{\kod}{\tfrac{k}{d}}
\newcommand{\nod}{\tfrac{n}{d}}
\newcommand{\ond}{\tfrac{1}{d}}
\newcommand{\twd}{\tfrac{2}{d}}
\newcommand{\thd}{\tfrac{3}{d}}
\newcommand{\onf}{\tfrac{1}{f}}
\newcommand{\Gr}{{\rm Gr}}
\newcommand{\Om}{\Omega}
\newcommand{\om}{\omega}

\newcommand{\Sp}{{\rm Sp}}
\newcommand{\bl}{\bigl}
\newcommand{\br}{\bigr}
\newcommand{\into}{\hookrightarrow}
\newcommand{\simto}{\buildrel\sim\over\longrightarrow}
\newcommand{\onto}{\mathop{\rlap{$\to$}\hskip2pt\h{$\to$}}}
\newcommand{\ges}{\geqslant}
\newcommand{\les}{\leqslant}
\newcommand{\bk}{\!\!\!}
\newcommand{\1}{\hskip1pt}
\newcommand{\mmu}{\rlap{$\mu$}\hskip1.1pt\h{$\mu$}}
\newcommand{\sw}{\,{\wedge}\,}
\newcommand{\ix}{\iota_{\xi}}
\begin{document}
\title[Roots of Bernstein-Sato polynomials]
{Roots of Bernstein-Sato polynomials of certain homogeneous polynomials with two-dimensional singular loci}
\author{Morihiko Saito}
\address{RIMS Kyoto University, Kyoto 606-8502 Japan}
\dedicatory{Dedicated to Professor Gert-Martin Greuel}
\begin{abstract} For a homogeneous polynomial of $n$ variables, we present a new method to compute the roots of Bernstein-Sato polynomial supported at the origin, assuming that general hyperplane sections of the associated projective hypersurface have at most weighted homogeneous isolated singularities. Calculating the dimensions of certain $E_r$-terms of the pole order spectral sequence for a given integer $r\in[2,n]$, we can detect its degeneration at $E_r$ for certain degrees. In the case of strongly free, locally positively weighted homogeneous divisors on ${\mathbb P}^3$, we can prove its degeneration almost at $E_2$ and completely at $E_3$ together with a symmetry of a modified pole-order spectrum for the $E_2$-term. These can be used to determine the roots of Bernstein-Sato polynomials supported at the origin, except for rather special cases.
\end{abstract}
\maketitle
\bs
\centerline{\bf Introduction}
\bsn
Let $f$ be a homogeneous polynomial of $n$ variables. Set $Z:=\{f=0\}\subset\PP^{n-1}$, and $d:=\deg f$. We assume $n,d\ges3$ in this paper. Let $b_f(s)$ be the Bernstein-Sato polynomial of $f$. We denote by $\R_f\subset\Q_{>0}$ the set of roots of $b_f(s)$ up to sign (see \cite{Ka}), since $R_f$ is confusing with a localization of $R=\C[x]$. For $z\in Z$, the local Bernstein-Sato polynomial $b_{h_z,z}(s)$ is {\it independent} of a choice of a local defining holomorphic function $h_z$ of $(Z,z)$, and is denoted by $b_{Z,z}(s)$, see for instance \cite{wh}. Let $\R_{Z,z}\subset\Q_{>0}$ be the set of roots of $b_{h_z,z}(s)$ up to a sign. Set
$$\R_Z=\mcup_{z\in Z}\,\R_{Z,z}\subset\R_f,$$
where the last inclusion holds by the above independence. Define
$$\R_f^0:=\R_f\setminus\R_Z\q\q\h{so that}\q\q\R_f=\R_f^0\sqcup\R_Z.$$
We call $\R_f^0$ the {\it roots of $b_f(s)$ up to a sign supported at the origin}.
\sk
Set $\Ff:=f^{-1}(1)\subset\C^n$ (the Milnor fiber). We have the monodromy eigenspaces of the Milnor cohomology
$$H^j(\Ff,\C)_{\la}:={\rm Ker}(T_s-\la)\q\h{for}\q\la\in\mmu_d:=\{\la\in\C\mid\la^d=1\},$$
where $T_s$ is the semisimple part of the monodromy. Note that $H^j(\Ff,\C)_{\la}=0$ for $\la\notin\mmu_d$. These spaces have the {\it pole order filtration} $P$, see \cite{Di1}, \cite{DiSa2}, etc. We have the following (see also Corollaries~1 and 2 below):
\msn
{\bf Theorem 1} (\cite[Theorem 2]{bfun}). {\it For $\al\notin\R_Z$, we have
$$\al\in\R_f^0\q\h{if}\q\Gr_P^p\,H^{n-1}(\Ff,\C)_{\ee(-\al)}\ne 0\q\bl(p=[n-\al]\br),
\leqno(1)$$
where $\ee(-\al):=e^{-2\pi i\al}$, and the converse holds if we have}
$$\al\notin\R_Z+\Z_{<0}.
\leqno(2)$$
\ms
In the notation explained below, the hypothesis in (1) is equivalent to that $\mu^{(\infty)}_k\ne 0$ for $\kod=\al$. Note that there are examples of {\it non-reduced\1} hyperplane arrangements of 3 variables where condition~(2) is {\it unsatisfied\1} and the converse of (1) {\it fails,} see Example 4.5 below.
\sk
In this paper we assume $\dim{\rm Sing}\,Z\les 1$, that is, $\dim{\rm Sing}\,f^{-1}(0)\les 2$, and moreover
$$\aligned&\raise-1mm\hbox{General hyperplane sections of $Z$ have at most}\\&\hbox{weighted homogeneous isolated singularities.}\endaligned
\leqno{\rm (GH)}$$
\sk
Let $\Om^{\ssb}$ be the graded complex of algebraic differential forms on $\C^n$. Its components are finite free graded modules over $R:=\C[x_1,\dots,x_n]$, where the $x_i$ are the coordinates of $\C^n$, which have degree 1 as well as the $\ddd x_i$, see (1.1) below. For $k\in\Z$, we have the microlocal {\it pole order spectral sequence}
$$E_1^{p,q}(f)_k=H^{p+q}_{\df\wedge}(\Om^{\ssb})_{qd+k}\Longrightarrow \Ht^{p+q-1}(\Ff,\C)_{\ee(-k/d)},
\leqno(3)$$
and its abutment filtration coincides with the pole order filtration $P$ in Theorem~1 up to the shift of filtration by $\bl[n-\kod\br]$, where $d:=\deg f$, see \cite{Di1}, \cite{DiSa1}, \cite{DiSa2} and (1.1.2) below. (Here the reader may assume $k\in[1,d]$ in (3) if he prefers.) Set
$$\aligned M:=H^n_{\df\wedge}(\Om^{\ssb}),&\q M^{(2)}:=H^n_{\ddd}(H^{\ssb}_{\df\wedge}(\Om^{\ssb})),\\
N:=H^{n-1}_{\df\wedge}(\Om^{\ssb})(-d),&\q N^{(2)}:=H^{n-1}_{\ddd}(H^{\ssb}_{\df\wedge}(\Om^{\ssb}))(-d),\\
Q:=H^{n-2}_{\df\wedge}(\Om^{\ssb})(-2d),&\q Q^{(2)}:=H^{n-2}_{\ddd}(H^{\ssb}_{\df\wedge}(\Om^{\ssb}))(-2d),\endaligned$$
where $M,N,Q$ are also denoted by $M^{(1)},N^{(1)},Q^{(1)}$ respectively. Recall that $(m)$ denotes the shift of grading by $m\in\Z$, that is, $G(m)_k=G_{k+m}$ ($k\in\Z$) for any graded module $G$.
\sk
For $r\ges 3$, we can define $M^{(r)}$, $N^{(r)}$, $Q^{(r)}$ using (1.1.4) below so that
$$E_r^{p,q}(f)_k=\begin{cases}
M^{(r)}_{qd+k}&\h{if}\,\,\,p+q=n,\\
N^{(r)}_{(q+1)d+k}\raise15pt\h{}&\h{if}\,\,\,p+q=n-1,\\
Q^{(r)}_{(q+2)d+k}\raise15pt\h{}&\h{if}\,\,\,p+q=n-2,\\
\end{cases}$$
where $E_r^{p,q}(f)_k=0$ otherwise (since $\dim{\rm Sing}\,f^{-1}(0)\les 2$). The formula holds also for $r=1,2$. It becomes quite simple if we restrict to $p=n$. The differentials
$$E_r^{n-1-q-r,\,q-1+r}(f)_k\buildrel{\ddd_r}\over\longrightarrow E_r^{n-1-q,\,q}(f)_k\buildrel{\ddd_r}\over\longrightarrow E_r^{n-1-q+r,\,q+1-r}(f)_k$$
are identified with
$$Q^{(r)}_{(q+1+r)d+k}\buildrel{{}''\!\ddd^{(r)}}\over\longrightarrow N^{(r)}_{(q+1)d+k}\buildrel{{}'\!\ddd^{(r)}}\over\longrightarrow M^{(r)}_{(q+1-r)d+k}\q\q(r\ges 1),
\leqno(4)$$
and $Q^{(r+1)}$, $N^{(r+1)}$, $M^{(r+1)}$ are identified respectively with
$${\rm Ker}\,{}''\!\ddd^{(r)},\q{\rm Ker}\,{}'\!\ddd^{(r)}/{\rm Im}\,{}''\!\ddd^{(r)},\q{\rm Coker}\,{}'\!\ddd^{(r)}\q\q(r\ges 1).$$
Set
$$\mu^{(r)}_k:=\dim M^{(r)}_k,\q\nu^{(r)}_k:=\dim N^{(r)}_k,\q\rho^{(r)}_k:=\dim Q^{(r)}_k\q\q(r\ges 1).$$
We will denote $\mu^{(1)}_k$, $\nu^{(1)}_k$, $\rho^{(1)}_k$ also by $\mu_k$, $\nu_k$, $\rho_k$ respectively (as in \cite{DiSa2}, \cite{wh}). Set $\mu_k^{(\infty)}:=\mu_k^{(r)}$ for $r\gg 0$ (note that $\{\mu_k^{(r)}\}$ is a weakly decreasing sequence for $r\ges 1$ with $k$ fixed), and similarly for $\nu_k^{(\infty)}$, $\rho_k^{(\infty)}$.
\sk
For $\,r\ges 2$, $\,\la\in\mmu_d\,$, $\,\beta\in\ond\1\Z_{>0}\,$, let $\,k\in[1,d]\,$ with $\,\la=\ee(-k/d)$, and define
$$\aligned&\chi^{(r)}_{f,\la}(\beta):=b^{(r),\beta}_{n-1,\la}-b^{(r),\beta}_{n-2,\la}+b^{(r),\beta}_{n-3,\la}+(-1)^{n-1}\delta_{\la,1}\q\q\q\h{with}\\
&b^{(r),\beta}_{n-1,\la}:=\sum_{j=0}^{[n-\beta-\kods]}\mu_{k+jd}^{(r)}\,,\,\,\,\,b^{(r),\beta}_{n-2,\la}:=\sum_{j=1}^{[n-\beta-\kods]}\nu_{k+jd}^{(r)}\,,\,\,\,\,b^{(r),\beta}_{n-3,\la}:=\sum_{j=2}^{[n-\beta-\kods]}\rho_{k+jd}^{(r)}\,,\endaligned$$
where $\delta_{\la,1}=1$ if $\la=1$, and $0$ otherwise. About the range of summation, note that
$$j\les\bl[n-\beta-\kod\br]\,\,\,{\Longleftrightarrow}\,\,\,k+jd\les d(n-\beta).$$
We call $\chi^{(r)}_{f,\la}(\beta)$, $b^{(r),\beta}_{n-1-i,\la}$ respectively the $(\la,r,\beta)$-{\it Euler characteristic} of the Milnor fiber and the $(\la,r,\beta)$-{\it Betti numbers} of the Milnor fiber ($i=0,1,2$).
\sk
The following is essential for a {\it certified\,} computation of the pole order spectral sequence.
\msn
{\bf Theorem~2.} {\it Assume condition~{\rm (GH)} holds, and there is an integer $r\in[2,n]$ together with $\beta\in\ond\1\Z_{>0}$ satisfying
$$\beta\les\bl\lceil\min\bl(\nod,\,\al_Z\br)\br\rceil^{/d}\,\,\,\,\h{with}\,\,\,\,\al_Z:=\min\R_Z,
\leqno(5)$$
$$\chi^{(r)}_{f,\la}\bl(\beta\br)\les(-1)^{n-1}\chi(U)\q\h{for}\,\,\,\,\la\in\mmu_d\,,
\leqno(6)$$
$$\mu^{(r)}_{k-md}\,\nu^{(r)}_k=0\q\h{if}\,\,\,\,\kod\les n-\beta,\,\,\,m\ges r\,,
\leqno(7)$$
where $\lceil\al\rceil^{/d}:=\min\bl\{\kod\in\ond\1\Z\mid\kod\ges\al\br\}$.
Then equality holds in $(6)$, and the pole order spectral sequence degenerates at $E_r$ for the highest cohomology and $\,\kod\les n-\beta$, that is,
$$\mu_k^{(r)}=\mu_k^{(\infty)}\q\h{for any}\,\,\,\,\kod\les n-\beta.
\leqno(8)$$
Moreover, condition~$(8)$ with $\mu_k$ replaced by $\nu_k$ and also by $\rho_k$ is valid for $\kod\les n-\beta$, if furthermore condition~$(7)$ with $\mu^{(r)}_{k-md}\,\nu^{(r)}_k$ replaced by $\nu^{(r)}_{k-md}\,\rho^{(r)}_k$ holds for $\kod\les n-\beta$.}
\ms
For an improvement of the assertions in certain cases, see Theorem~(4.3) below. Assuming condition~(GH) (and using \cite{wh}), the assumption in the last assertion of Theorem~2 may be replaced by the following:
$$\h{$Z$ is locally analytically trivial along general points of ${\rm Sing}\,Z$,}
\leqno{\rm (AT)}$$
see also Remark~(2.5)(i) below.
\sk
In (5) we have
$$\min\bl(\nod,\,\al_Z\br)=\al_f\,(:=\min\R_f),$$
see (2.1.1) below. By \cite[Corollary~3.6]{dFEM} (see Remark~(1.8) below), we have
$$\bl\lceil\min\bl(\nod,\,\al_Z\br)\br\rceil^{/d}\ges\thd.$$
So condition~(5) is satisfied if $\beta\les\thd$. The assertion (8) for $\mu_k$ is sufficient to get the roots of Bernstein-Sato polynomials using Theorem~1. Condition~(7) trivially holds if $r=n$. (Recall that the differentials ${}'\!\ddd^{(r)}$, ${}''\!\ddd^{(r)}$ shift the degree by $-rd$, see (4).) Set
$$\Lambda(r):=\bl\{\ee(-k/d)\,\,\big|\,\,\mu^{(r)}_{k-md}\,\nu^{(r)}_k\ne0\,\,\,\bl(\kod\les n-\beta,\,\exists\,m\ges r\br)\br\}\subset\mmu_d.$$
Condition~(7) can be replaced in Theorem~2 by the following.
$$b^{(r),\beta}_{n-1,\la}=b^{(r),\beta}_{n-1,\overline{\la}}\,\,\,\,(\forall\,\la\in\mmu_d),\q\h{and}\q\Lambda(r)\cap\overline{\Lambda(r)}=\emptyset.
\leqno(9)$$
The first condition of (9) is {\it not\1} automatically satisfied (since the relation between the pole order spectral sequence and the {\it real structure\1} of the Milnor cohomology is unclear) although it should hold if the spectral sequence degenerates at $E_r$ for the highest cohomology and $\kod\les n-\beta$. The second condition of (9) is satisfied if $d$ is rather small.
It seems quite difficult to show the equalities $\mu_k^{(r)}=\mu_k^{(\infty)}$ for $\kod\les n-\beta$ without calculating also the $\nu_k^{(r)}$, $\rho_k^{(r)}$ (compare \cite{DiSt2}).
\sk
For the proof of Theorem~2, we need the following.
\msn
{\bf Proposition~1.} {\it Assume condition~{\rm (GH)} holds. Let $\alt'_Z:=\min\,\Rt_{Z'}$ in the notation of $(2.1)$ below with $Z'$ a general hyperplane section of $Z$. Then}
$$\aligned&{}''\!\ddd^{(1)}:Q_k\to N_{k-d}\,\,\,\,\h{\it is injective and}\,\,\,\,Q^{(r)}_k=0\,\,\,(\forall\,r\ges 2),\\&\q\h{\it if}\,\,\,\,\kod>n-\lceil\alt'_Z\rceil^{/d}\,.\endaligned
\leqno(10)$$
\sk
Note that the assumption in (10) is equivalent to $\kod>n-\alt'_Z$, and we have
$$\lceil\alt'_Z\rceil^{/d}\ges\thd,$$
by \cite[Corollary~3.6]{dFEM}, see also Remark~(1.8) below. Proposition~1 is proved by using \cite[Theorem~5.3]{DiSa2} (applied to general hyperplane sections of $Z$) together with a kind of ``torsion-freeness" of $Q$, see (3.4.1) and (2.3) below. 
\sk
We now restrict to a special case including reduced hyperplane arrangements and strongly free, locally positively weighted homogeneous divisors (see (4.1) below) on $\PP^3$. For $r\ges 2$, $\la\in\mmu_d$, $\beta\in\ond\1\Z_{>0}$, define the {\it modified $(\la,r,\beta)$-Euler characteristic} of the Milnor fiber by
$$\aligned&\cho^{(r)}_{f,\la}(\beta):=\bo^{(r),\beta}_{3,\la}-\bo^{(r),\beta}_{2,\la}+b^{(r),\beta}_{1,\la}-\delta_{\la,1}\q\q\q\q\q\h{with}\\
&\bo^{(r),\beta}_{3,\la}:=\sum_{j=0}^{[2-\beta-\kods]}\mu_{k+jd}^{(r)}\,,\q \bo^{(r),\beta}_{2,\la}:=\sum_{j=1}^{[3-\beta-\kods]}\nu_{k+jd}^{(r)}\,,\q b^{(r),\beta}_{1,\la}:=\sum_{j=2}^{[4-\beta-\kods]}\rho_{k+jd}^{(r)}\,,\endaligned$$
where $k\in[1,d]$ with $\la=\ee(-k/d)$ as above. Note that the definition of $b^{(r),\beta}_{1,\la}$ is the same as in Theorem~2 with $n=4$. We have a variant of Theorem~2 as follows.
\msn
{\bf Theorem~3.} {\it Assume $n=4$, condition~{\rm (GH)} holds, and there is an integer $r\in[2,4]$ with following conditions satisfied\,$:$
$$\max\R_f\les 2-\twd\,,
\leqno(11)$$
$$\cho^{(r)}_{f,\la}\bl(\twd\br)\les-\chi(U)\q\h{for}\,\,\,\,\la\in\mmu_d\,,
\leqno(12)$$
$$\mu^{(r)}_{k-md}\,\nu^{(r)}_k=0\q\h{if}\,\,\,\,\kod\les 3-\twd,\,\,\,m\ges r\,,
\leqno(13)$$
where $\beta=\twd$. Then equality holds in $(12)$, the pole order spectral sequence degenerates at $E_r$ for the highest cohomology with $\kod\les 2-\twd$, that is, $(8)$ holds for $\kod\les 2-\twd$, and we have
$$\nu^{(r)}_k=0\q\h{for}\q\kod\in\bl(3-\twd,4-\twd\br].
\leqno(14)$$
Moreover an analogue of the last assertion of Theorem~$2$ also holds, that is, if we assume condition~$(13)$ with $\mu^{(r)}_{k-md}\,\nu^{(r)}_k$ replaced by $\nu^{(r)}_{k-md}\,\rho^{(r)}_k $ for $\kod\les 4-\twd$, $m\ges r$, then we get the assertion~$(8)$ with $\mu_k$ replaced by $\nu_k$, and also by $\rho_k$, for $\kod\les 4-\twd$.}
\ms
In the case of strongly free, locally positively weighted homogeneous divisors on $\PP^3$, we can show the degeneration of the pole order spectral sequence almost at $E_2$ and completely at $E_3$ together with a symmetry of the modified pole-order spectrum for the $E_2$-term, see Theorems~(4.5) and (4.9) below.
Condition~(11) implies by Theorem~1 above that
$$\mu^{(\infty)}_k=0\q\h{if}\q\kod>2-\twd.
\leqno(15)$$
This condition is satisfied if $Z$ is an hyperplane arrangement or a strongly free, locally weighted homogeneous divisor on $\PP^{n-1}$, see Remark~(4.4)(ii) below. For these divisors with $n=4$, it seems to be expected that the pole order spectral sequence would degenerate at $E_2$, see \cite{DiSt2}. It seems, however, quite possible that the hypotheses of Theorem~3 are satisfied for $r=2$ except for condition~(7) modified as in Theorem~3, while the latter assumption is satisfied only for $r=3$ so that we can get the partial $E_2$-degeneration only after calculating certain $E_3$-terms. (In some simple cases as Example~(5.8) below, all the hypotheses of Theorem~3 are satisfied already for $r=2$.)
\sk
Theorems~2--3 and (4.10) below may be used to show that the results of some computations in \cite{DiSt2} are correct without assuming a conjecture, for instance, in the case of Example~5.7 in {\it loc.\,cit.}, see Example~(5.8) below.
\sk
In this paper we also present a new algorithm to compute the $\mu^{(r)}_k$, $\nu^{(r)}_k$, $\rho^{(r)}_k$, see (1.2) below. In the case $\dim{\rm Sing}\,Z=0$, the $\rho^{(r)}_k$ vanish and the $\nu^{(r)}_k$ can be determined by the $\mu^{(r)}_k$ (see for instance \cite{DiSa2}) so that the calculation is much easier than the case $\dim{\rm Sing}\,Z=1$. Our algorithm is more systematic than the one in \cite{DiSt1} for the case $\dim{\rm Sing}\,Z=0$, where ``syzygies" are used to calculate the {\it kernel\,} of the differential $\ddd_2$ rather than its {\it image} as in our paper (and a simple version of Theorem~2 is actually hired there to make the calculation {\it reliable}). In \cite{DiSt2}, ``syzygies" are used in a more intelligent way than in \cite{DiSt1} so that the computing time of some $E_2$-terms is reduced in certain cases with $\dim{\rm Sing}\,Z=1$. (This could be applied also to the computation of certain $E_3$-terms in the case where $\dim{\rm Sing}\,Z=0$ and $n=3$, see Remark~(1.7)(iii) below.) However, a major problem seems to be a lack of methods to assure the validity of all the computations in \cite{DiSt2}, see Remark~(1.7)(ii) below.
\sk
It seems rather complicated to employ ``syzygies" also for reducing the computing time of $\nu^{(r)}_k$, $\rho^{(r)}_k$ in Theorems~2 and 3 for the case $\dim{\rm Sing}\,Z=1$, $n=4$. It may take more than several hours to compute the $\mu^{(r)}_k$, $\nu^{(r)}_k$, $\rho^{(r)}_k$ ($r\les 3$) even for the case $n=4$, $d=5$ in general (except for the hyperplane arrangement case), where the actual computing time may depend on each example since only the terms with condition~(7) unsatisfied for $m<r$ are calculated. It seems, however, more important to find a sufficient condition for a partial degeneration as in Theorem~2 in order to make all the calculations {\it certified\,} in \cite{DiSt2}.
\sk
Combining Theorems~1 and 2 (or 3), we can determine $\R_f^0$ in many cases, although there is a difficulty coming form the assumption~(2) in Theorem~1. To avoid this problem, set as in \cite{wh}
$${\rm CS}(f):=\bl\{k\in\Z\,\,\big|\,\,\kod\in[\al_Z,\,n-2-\al_Z]\cap(\R_Z+\Z_{<0})\setminus\R_Z\br\}\subset\Z.$$
Here we use the inclusion $\Rt_Z\subset[\alt_Z,n-1-\alt_Z]$ in the notation of (2.1) below, see \cite[Theorem~0.4]{mic} (and (2.1.2) below).
\sk
By an argument similar to the proof of \cite[Theorem 4]{wh} (with $\delta_k:=\mu_k-\nu_{k+d}$ replaced by $\mu^{(\infty)}_k\,$), Theorem~1 implies the following.
\msn
{\bf Corollary~1.} {\it Assume $\,\min{\rm CS}(f)\ges n$, and moreover
$$\mu^{(\infty)}_k>0\,\,\,\,\,\h{for any}\,\,\,\,\,k\in{\rm CS}(f).
\leqno(16)$$
Then $\R^0_f$ can be determined by the $\mu^{(\infty)}_k$ {$(k\in\N)$}. More precisely, we have the inclusion
$$\R_f^0\subset\bl[\nod,n\br)\cap\ond\,\Z,
\leqno(17)$$
and for $\kod\in\bl(\bl[\nod,n\br)\cap\ond\,\Z\br)\setminus\R_Z$, we have $\kod\in\R^0_f$ if and only if $\,\mu^{(\infty)}_k>0$.}
\ms
The first assumption follows from condition~(16), since $\mu_k^{(\infty)}=0$ for $k<n$. It is known that condition~(16) does not hold for certain extremely degenerated curves having only weighted homogeneous singularities, see \cite{wh}. There are also examples with condition~(16) unsatisfied in the case $Z$ has non-quasi-homogeneous isolated singularities, see Example~(5.3) below. For the moment such examples are restricted to the case $\chi(U)\les 1$ if $n=3$.
\msn
{\bf Remark~1.} The method in this paper is {\it not very efficient} in the case $Z$ has only {\it weighted homogeneous isolated singularities,} since the $E_2$-degeneration and a certain injectivity of the differential $\ddd_1$ hold by \cite[Theorem 2]{wh} and \cite[Theorem~5.3]{DiSa2} so that it is enough to compute the {\it Hilbert series of the graded Milnor algebra}, which gives the $E_1$-term of the spectral sequence, see \cite{wh}.
\msn
{\bf Remark~2.} In the non-quasi-homogeneous isolated singularity case, it is quite nontrivial to determine $\R_Z$ by calculating the local Bernstein-Sato polynomials of $Z$, where a computer program like RISA/ASIR (see \cite{NoTa}) is needed. This is quite different from the {\it weighted homogeneous isolated singularity} case, where the local Bernstein-Sato polynomials can be determined only by the {\it weights} as a consequence of results of \cite{Ma}, \cite{Sat}, \cite{ScSt}, \cite{Va1}, see \cite[Section 1.9]{wh}.
\ms
Set ${\rm Supp}\,\{\mu^{(\infty)}_k\}:=\{k\in\Z\mid\mu^{(\infty)}_k\ne 0\}$, which is called the {\it support} of the $\mu^{(\infty)}_k$. We say that the support of the $\mu^{(\infty)}_k$ is {\it discretely connected outside} $\R_Z$ if the following holds:
$$\aligned&{\rm Supp}\,\{\mu^{(\infty)}_k\}\setminus d\,\R_Z=\bl([k_{\rm min},k_{\rm max}]\cap\Z\br)\setminus d\,\R_Z,\\ &\h{with}\q\q
k_{\rm min},\,k_{\rm max}\in{\rm Supp}\,\{\mu^{(\infty)}_k\}\setminus d\,\R_Z.\endaligned
\leqno(18)$$
We do not know any example with this condition unsatisfied. There is an analogue of \cite[Corollary~3]{wh} (which is also a corollary of Theorem~1) as follows:
\msn
{\bf Corollary~2.} {\it Assume ${\rm CS}(f)\subset[n,\infty)$, condition $(18)$ holds, and moreover
$$\bl([k_{\rm min},k_{\rm max}]\cap\Z\br)\setminus d\,\R_Z=\bl([n,k_{\rm max}]\cap\Z\br)\setminus d\,\R_Z,
\leqno(19)$$
$$k_{\rm max}\ges\max\bl(d\,\R_Z\cap\Z\br)-d.
\leqno(20)$$
Then we have}
$$d\,\R_f^0=\bl([n,k_{\rm max}]\cap\Z\br)\setminus d\,\R_Z.$$
\ms
This does not necessarily imply that $d\,\R_f\cap\Z=[n,k']\cap\Z$ for some $k'\in\N$ (since it might be something like $\bl([n,k']\cap\Z\br)\cup\{k'+2\}$ with $k'+2\in d\,\R_Z$). Condition~(20) is needed to show that $\kod\notin\R_f^0$ for $k>k_{\rm max}$ using Theorem~1. We can replace it by the condition that $k_{\rm max}\ges d-1$, if $n=3$. Conditions~(18) and (20) are always satisfied as far as calculated (with $n=3$).
\sk
As for condition (19), the inequality $k_{\rm min}\ges n$ always holds, and the strict inequality occurs {\it only if} $\al_Z\les\nod$, see \cite[4.2.6]{wh}. In this case we {\it cannot\1} determine whether $\kod\in\R_f^0$ for $k\in[n,k_{\rm min}-1]\cap\Z$ with $\kod\in\bl(\R_Z+\Z_{<0}\br)\setminus\R_Z$ because of condition~(2) in Theorem~1.
\sk
This work is partially supported by Kakenhi 15K04816. We thank A.~Dimca for useful discussions about the subject of this paper, and also for making a computer program based on the algorithm in this paper (which encouraged us to make a less sophisticated one). We thank G.-M.~Greuel for useful comments related to the codes of Singular in this paper.
Finally we thank the referees for correcting some errors and improving the paper.
\sk
In Section~1 we review some basics of pole order spectral sequences, and present a new algorithm for computations of the spectral sequence. In Section~2 we prove Theorems~2 and 3 together with Propositions~1. In Section~3 we prove a structure theorem together with self-duality isomorphisms for the Koszul cohomologies $M,N,Q$. These imply a kind of ``torsion-freeness" of $Q$. In Section~4 we prove Theorems~(4.3--4) and (4.10) improving Theorems~2--3 in certain cases, and show a symmetry of a modified pole-order spectrum for strongly free divisors in Theorem~(4.5). In Section~5 we calculate some examples.
In Appendix we show a {\it double symmetry\1} of the modified pole-order spectrum for the $E_2$-term in the case of strongly free divisors on $\PP^3$.
\bs\bs
\vbox{\centerline{\bf 1. Pole order spectral sequences}
\bsn
In this section we review some basics of pole order spectral sequences, and present a new algorithm for computations of the spectral sequence.}
\msn
{\bf 1.1.~Spectral sequence.} Let $R=\C[x_1,\dots,x_n]$, and $\Om^j:=\Gamma(\C^n,\Om_{\C^n}^j)$ ($j\in\Z$), where $x_1,\dots,x_n$ are the coordinates of $\C^n$, and the $\Om_{\C^n}^j$ are algebraic so that the $\Om^j$ are finite free graded $R$-modules with $\deg x_i=\deg\ddd x_i=1$. For instance,
$$\Om^n=R\,\ddd x_1\wedge\cdots\wedge\ddd x_n.$$
In the notation of the introduction, the pole order spectral sequence is essentially the spectral sequence associated with the graded double complex $C^{\ssb,\ssb}$ with
$$C^{p,q}=\Om^{p+q}(qd),
\leqno(1.1.1)$$
where $(qd)$ denotes the shift of grading, that is, $\Om^{p+q}(qd)_k=\Om^{p+q}_{k+qd}$. The anti-commuting two differentials are given by
$$\ddd:C^{p,q}\to C^{p+1,q},\q\df\wedge:C^{p,q}\to C^{p,q+1}\q\q(p,q\in\Z),$$
where the first differential $\ddd$ is the usual one.
\sk
We have the decreasing filtration $F$ on the total complex $C^{\ssb}$ defined by
$$F^pC^j=\mopl_{i\ges p}\,C^{i,j-i}.$$
(In this section, the pole order filtration $P$ is denoted by $F$.)
This gives the microlocal {\it pole order spectral sequence} (see \cite{Di1}, \cite{DiSa2}):
$$E_1^{p,q}(f)_k=H^{p+q}_{\df\wedge}(\Om^{\ssb})_{qd+k}\Longrightarrow H^{p+q}C^{\ssb}_k=\Ht^{p+q-1}(\Ff,\C)_{\ee(-k/d)},
\leqno(1.1.2)$$
where $\ee(-\al)$ is as in Theorem~1, and the last isomorphism of (1.1.2) follows from the graded version of the theory of microlocal Gauss-Manin systems (see \cite{BaSa} for the analytic case). Its abutment filtration coincides with the pole order filtration $P$ in Theorem~1 up to the shift of filtration by $\bl[n-\kod\br]$, see \cite{Di1} (and \cite[Section 1.8]{DiSa1} for the case of the top cohomology). For the relation with the {\it usual\,} pole order spectral sequence, see \cite{DiSa2}.
\sk
For $m\in\Z$, there are filtered graded isomorphisms
$$(C^{\ssb},F)\simto(C^{\ssb}(md),F[m]),
\leqno(1.1.3)$$
induced by the natural isomorphisms
$$C^{p,q}_k\simto C^{p+m,q-m}_{k+md}\q(p,q,k\in\Z).$$
(Indeed, the condition $i\ges p$ implies $i':=i+m\ges p+m$, and $(F[m])^p=F^{p+m}$.)
\sk
These imply the isomorphisms
$$E_r^{p,q}(f)_k\simto E_r^{p+m,q-m}(f)_{k+md}\q(p,q,k\in\Z,\,r\ges 1).
\leqno(1.1.4)$$
\msn
{\bf 1.2.~Computation of the differential $\ddd_r$.} By the general theory of spectral sequences together with the isomorphisms (1.1.3), we have a decreasing (non-separated) filtration $Z^{\ssb}$ and an increasing (non-exhaustive) filtration $B_{\ssb}$ on
$$C^{\,n,-i}=\Om^{n-i}(-id)\q\q(i\in[0,n]),$$
such that
$$B_{r-1}\Om^{n-i}\subset B_r\Om^{n-i}\subset B_{\infty}\Om^{n-i}\subset Z^{\infty}\Om^{n-i}\subset Z^{r+1}\Om^{n-i}\subset Z^r\Om^{n-i},
\leqno(1.2.1)$$
$$E_r^{n,-i}(f)=(Z^r\Om^{n-i}/B_{r-1}\Om^{n-i})(-id)\q\q(r\ges 1).
\leqno(1.2.2)$$
Moreover, using the isomorphisms (1.1.4), the differential $\ddd_r$ can be identified with the composition of an isomorphism
$$(\Gr_Z^r\Om^{n-i-1})(-(i+1)d+rd)\simto(\Gr_r^B\Om^{n-i})(-id)\q\q(r\ges 1),
\leqno(1.2.3)$$
with the surjection
$$Z^r\Om^{n-i-1}/B_{r-1}\Om^{n-i-1}\onto\Gr_Z^r\Om^{n-i-1},$$
and the injection
$$\Gr_r^B\Om^{n-i}\into Z^r\Om^{n-i}/B_{r-1}\Om^{n-i},$$
up to the shift of grading. So the coimage and the image of $\ddd_r$ are identified respectively with the source and the target of (1.2.3). Here
$$\aligned B_0\Om^{n-i}={\rm Im}\,\df\wedge,&\q Z^1\Om^{n-i}={\rm Ker}\,\df\wedge,\\
B_{-1}\Om^{n-i}=0,&\q Z^0\Om^{n-i}=\Om^{n-i}.\endaligned$$
Then it is enough to calculate $\dim\Gr^B_{\ssb}\Om^{n-i}_{k-id}$ for $0\les i<\dim{\rm Sing}\,f^{-1}(0)$ in order to know that of $\Gr^r_Z\Om^{n-i-1}_{k-(i+1)d+rd}$. (For $i=\dim{\rm Sing}\,f^{-1}(0)$, it is enough to calculate $B_0\Om^{n-i}$, since $\Gr^B_r\Om^{n-i}=0$ for $r>1$.)
\sk
We have the isomorphisms
$$\begin{array}{cccc}Z^r\Om^{n-i}_{k-id}/B_{r-1}\,\Om^{n-i}_{k-id}&=&E_r^{n,-i}(f)_k,\\ \cup&&\cup\,\,\\ \Gr^B_r\Om^{n-i}_{k-id}&=&{\rm Im}\,\ddd_r&(r\ges 1),\end{array}
\leqno(1.2.4)$$
and the filtration $B_{\ssb}$ can be described explicitly as follows.
\sk
For $r\ges 0$, $k\in\Z$, we have the morphism

$$\aligned&\begin{array}{rcl}\Psi^{(r)}_{i,k}\,:\,\mopl_{j=-1}^{r-1}\,\Om^{n-i-1}_{jd+k-id}&\to&\mopl_{j=0}^r\,\Om^{n-i}_{jd+k-id}\\
\rlap{$\hskip.26pt{\scriptstyle\,|}$}{\scriptstyle\cup}\q\q\q& &\q\q\rlap{$\hskip.26pt{\scriptstyle\,|}$}{\scriptstyle\cup}\\
(\om_{-1},\dots,\om_{r-1})&\mapsto&(\om'_0,\dots,\om'_r)\end{array}\\
\h{with}\raise8mm\h{}\q&\,\,\,\,\om'_j:=\begin{cases}\df{\wedge}\,\om_{j-1}+\ddd\om_j&\h{if}\,\,\,j=0,\dots,r-1,\\ \df{\wedge}\,\om_{j-1}&\h{if}\,\,\,j=r.\end{cases}\q\q\h{}\endaligned
\leqno(1.2.5)$$
\skn
Define
$$\Phi^{(r)}_{i,k}:\mopl_{j=-1}^{r-1}\,\Om^{n-i-1}_{jd+k-id}\to\mopl_{j=1}^r\,\Om^{n-i}_{jd+k-id}
\leqno(1.2.6)$$
to be the composition of $\Psi^{(r)}_{i,k}$ with the natural projection
$$\mopl_{j=0}^r\,\Om^{n-i}_{jd+k-id}\onto\mopl_{j=1}^r\,\Om^{n-i}_{jd+k-id},$$
(where $\Phi^{(0)}_{i,k}=0$). We have the following.
\msn
{\bf 1.3.~Proposition.} {\it In the above notation, there are canonical isomorphisms}
$$\Psi^{(r)}_{i,k}({\rm Ker}\,\Phi^{(r)}_{i,k})=B_r\,\Om^{n-i}_{k-id}\q\q\bl(r\ges 0,\,k\in\Z\br).
\leqno(1.3.1)$$
\msn
{\it Proof.} We may assume $r>0$, since the assertion easily follows from definitions if $r=0$.
By the definition of the spectral sequence associated with the filtered complex $(C^{\ssb},F)$, we can identify ${\rm Ker}\,\Phi^{(r)}_{i,k}$ with the quotient of
$$Z_{r,k}^{n-r,\,r-i-1}:={\rm Ker}\bl(\ddd:F^{n-r}C^{n-i-1}_k\to C^{n-i}/F^nC^{n-i}_k\br),$$
divided by
$$F^{n+1}C^{n-i-1}_k=\mopl_{j'<-i-1}\,C_k^{n-i-1-j',j'}=\mopl_{j<-1}\,\Om^{n-i-1}_{jd+k-id},$$
where $j'=j-i$, and we have the isomorphisms
$$\aligned F^{n-r}C^{n-i-1}_k&=\mopl_{j'\les r-i-1}\,C_k^{n-i-1-j',j'}=\mopl_{j\les r-1}\,\Om^{n-i-1}_{jd+k-id},\\ F^nC^{n-i}_k&=\mopl_{j'\les-i}\,C_k^{n-i-j',j'}=\mopl_{j\les 0}\,\Om^{n-i}_{jd+k-id}.\endaligned$$
Moreover the restriction of $\Psi^{(r)}_{i,k}$ to ${\rm Ker}\,\Phi^{(r)}_{i,k}$ is induced by the differential of the quotient complex $C^{\ssb}_k/F^{n+1}C^{\ssb}_k$. So the assertion follows from Remark below.
\msn
{\bf Remark.} By the general theory of spectral sequences associated with filtered complexes $(C^{\ssb},F)$ (see for instance \cite[Section 1.3]{De}, \cite{Go}), it is well known that $E^{p,q}_r$ is a quotient of
$$Z^{p,q}_r:={\rm Ker}(\ddd:F^pC^{p+q}\to C^{p+q+1}/F^{p+r}C^{p+q+1}),$$
and the differential $\ddd_r$ is induced by the restriction of the differential $\ddd$ of $C^{\ssb}$.
\msn
{\bf 1.4.~Corollary.} {\it In the notation of $(1.2)$, set
$$\beta^{(r)}_{i,k}:={\rm rank}\,\Psi^{(r)}_{i,k}-{\rm rank}\,\Phi^{(r)}_{i,k}\q\q(r\ges 0).
\leqno(1.4.1)$$
Then, for $i=0$, we have the following equality in the notation of the introduction}
$$\mu_k^{(r+1)}=\tbinom{k-1}{n-1}-\beta^{(r)}_{0,k}\q\q(r\ges 0).
\leqno(1.4.2)$$
\msn
{\it Proof.} For $r\ges 0$, we have $\dim \Om^n_k=\tbinom{k-1}{n-1}$, $Z^r\Om^n=\Om^n$ (since $\Om^{n+1}=0$), and
$$\aligned\dim\Psi^{(r)}_{i,k}({\rm Ker}\,\Phi^{(r)}_{i,k})&=\dim{\rm Ker}\,\Phi^{(r)}_{i,k}-\dim{\rm Ker}\,\Psi^{(r)}_{i,k}\\&={\rm rank}\,\Psi^{(r)}_{i,k}-{\rm rank}\,\Phi^{(r)}_{i,k}=\beta^{(r)}_{i,k}.\endaligned
\leqno(1.4.3)$$
Here the first equality follows from the inclusion ${\rm Ker}\,\Psi^{(r)}_{i,k}\subset{\rm Ker}\,\Phi^{(r)}_{i,k}$, since this implies that
$${\rm Ker}\,\bl(\Psi^{(r)}_{i,k}\big|{\rm Ker}\,\Phi^{(r)}_{i,k}\br)={\rm Ker}\,\Psi^{(r)}_{i,k}.
\leqno(1.4.4)$$
So the assertion follows from Proposition~(1.3) together with (1.2.2).
\msn
{\bf 1.5.~Corollary.} {\it In the notation of the introduction and $(1.4.1)$, assume $\dim{\rm Sing}\,Z=0$. Then we have
$$\nu^{(r)}_k-\nu^{(r+1)}_k=\mu^{(r)}_{k-rd}-\mu^{(r+1)}_{k-rd}=\beta^{(r)}_{0,k-rd}-\beta^{(r-1)}_{0,k-rd}\q\q(r\ges 1).
\leqno(1.5.1)$$
Moreover, for $r=1$, we have}
$$\mu_k-\nu_k=\gamma_k\,\,\,(k\in\Z)\q\h{with}\q\msum_k\,\gamma_k\,v^k=\bl(v+\cdots+v^{d-1}\br)^n.
\leqno(1.5.2)$$
\msn
{\it Proof.} Considering the rank of $\ddd^{(r)}$, we get the first equality of (1.5.1), since $\Gr_r^B\Om^{n-1}=0$ ($r\ges 1$). The second equality of (1.5.1) follows from Corollary~(1.4). For (1.5.2), see for instance \cite[Formula (3)]{DiSa2}. This finishes the proof of Corollary~(1.5).
\msn
{\bf 1.6.~Corollary.} {\it In the notation of the introduction and $(1.4.1)$, assume $\dim{\rm Sing}\,Z=1$. Then we have the equalities}
$$\rho^{(r)}_k-\rho^{(r+1)}_k=\beta^{(r)}_{1,k-rd}-\beta^{(r-1)}_{1,k-rd}\q\q(r\ges 1).
\leqno(1.6.1)$$
$$\aligned\nu^{(r)}_k-\nu^{(r+1)}_k&=\mu^{(r)}_{k-rd}-\mu^{(r+1)}_{k-rd}+\rho^{(r)}_{k+rd}-\rho^{(r+1)}_{k+rd}\\&=\beta^{(r)}_{0,k-rd}-\beta^{(r-1)}_{0,k-rd}+\beta^{(r)}_{1,k}-\beta^{(r-1)}_{1,k}.\endaligned
\leqno(1.6.2)$$
\msn
{\it Proof.} These assertions follow from Proposition~(1.3) together with (1.2.2) and (1.4.3) by considering the ranks of the differentials $''\!\ddd^{(r)}$, $'\!\ddd^{(r)}$ in (4) in the introduction.
\msn
{\bf 1.7.~Remarks.} (i) In \cite[Section 4.1]{DiSt1} ``syzygies" are used to calculate the ``kernel" of the differential $\ddd_r$ ($r=1,2$), giving the $E_{r+1}$-term of the pole order spectral sequence for any reduced polynomial of 3 variables. For $r=1$, the method in our paper studying the image of $\ddd_1$ without using syzygies seems to be almost as fast as the one in \cite{DiSt1}, since the matrix whose rank we have to calculate (which is represented by $\phi'_q$ in \cite[Section~4.1]{DiSt1}) is essentially the same as far as the differential $\ddd_1$ is concerned.
\sk
For $r=2$, the method in our paper seems slightly better, since their method in \cite[Section~5.4]{DiSt1} has to compute the rank of the matrix corresponding to the linear map
$$\aligned&\Om^2_q\oplus \Om^2_{q-d}\oplus\Om^2_{q-d}\oplus\Om^2_{q-2d}\q\to\q\Om^3_{q+d}\oplus\Om^3_q\oplus\Om^3_q\oplus\Om^3_{q-d}\\
&\q\q\q\q\rlap{$\hskip.26pt{\scriptstyle\,|}$}{\scriptstyle\cup}\q\q\q\q\q\q\q\q\q\q\q\q\q\q\,\,\rlap{$\hskip.26pt{\scriptstyle\,|}$}{\scriptstyle\cup}\\
&\q(\om_0,\,\om_1,\,\om'_1,\,\om_2)\,\mapsto\,(\df{\wedge}\1\om_0,\,\ddd\om_0{-}\1\df{\wedge}\1\om_1,\,\df{\wedge}\1\om'_1,\,\ddd\om_1{-}\ddd\om'_1{-}\1\df{\wedge}\1\om_2)\endaligned
\leqno(1.7.1)$$
whereas it is enough to calculate the following in our algorithm (see (1.2.5)):
$$\aligned&\Om^2_q\oplus\Om^2_{q-d}\oplus\Om^2_{q-2d}\q\to\q\Om^3_{q+d}\oplus\Om^3_q\oplus\Om^3_{q-d}\\
&\q\q\q\q\rlap{$\hskip.26pt{\scriptstyle\,|}$}{\scriptstyle\cup}\q\q\q\q\q\q\q\q\q\q\q\rlap{$\hskip.26pt{\scriptstyle\,|}$}{\scriptstyle\cup}\\
&\q\q(\om_0,\,\om_1,\,\om_2)\q\mapsto\q(\df{\wedge}\1\om_0,\,\ddd\om_0\,{+}\,\df{\wedge}\1\om_1,\,\ddd\om_1\,{+}\,\df{\wedge}\1\om_2).\endaligned
\leqno(1.7.2)$$
It is possible to capture the rank of $\Phi^{(r)}_{i,k}$ in (1.2.6) during the computation of the rank of $\Psi^{(r)}_{i,k}$ in (1.2.5) provided that the program is well designed, since these correspond to matrices of size $(a',b)$ and $(a,b)$ with $a'\les a$. Hence the complexity of the calculation of the former can be neglected.
\ms
(ii) In \cite{DiSt2}, ``syzygies" are used quite efficiently for computations of the pole order filtration on the {\it highest\1} Milnor cohomology of certain homogeneous polynomials such that the $E_2$-degeneration of the pole order spectral sequences is expected, more precisely, if $f$ defines a hyperplane arrangement or a strongly free, locally positively weighted homogeneous divisor. Some of their computations, however, entirely depend on \cite[Conjecture 3.8]{DiSt2} which does not seem easy to prove in general.
\sk
This conjecture is proved in \cite[Theorem 3.9]{DiSt2} only for the {\it unipotent monodromy} part by using the theory of logarithmic forms, since the pole order filtration is {\it trivial} on this part so that the logarithmic forms are enough in the above two cases. However, this triviality does not necessarily hold on the {\it non-unipotent monodromy} part, and it seems quite difficult to figure out how to extend the above argument to the case where the pole order filtration is non-trivial and hence the logarithmic differential forms would not be enough.
\ms
(iii) It is easy to extend the new method in Remark~(ii) above to the calculation of the image of the differential $\ddd_2$ giving the $E_3$-term in the 3 variable case.
Indeed, if the ``syzygies" (that is, the {\it kernel}\1) of $\df\wedge:\Om^2\to\Om^3$ have generators $\eta_1,\dots,\eta_r$ over $R$ with pure degrees $d_1,\dots,d_r$, then it is enough to consider the following (instead of (1.7.2)):
$$\aligned&\bl(\mopl_{i=1}^r\,R_{q-d_i}\br)\oplus\Om^2_{q-d}\oplus\Om^2_{q-2d}\q\to\q\Om^3_q\oplus\Om^3_{q-d}\\
&\q\q\q\q\q\rlap{$\hskip.26pt{\scriptstyle\,|}$}{\scriptstyle\cup}\q\q\q\q\q\q\q\q\q\q\q\q\,\,\rlap{$\hskip.26pt{\scriptstyle\,|}$}{\scriptstyle\cup}\\
&\q(g_1\dots,g_r,\,\om_1,\,\om_2)\,\mapsto\,\bl(\msum_{i=1}^r\ddd(g_i\eta_i)\,{+}\,\df{\wedge}\1\om_1,\,\ddd\om\,{+}\,\df{\wedge}\1\om_2\br).\endaligned
\leqno(1.7.3)$$
Note that the first direct factor $\Om^3_{q+d}$ in the target of (1.7.2), which has the biggest dimension among the direct factors, is eliminated in the target of (1.7.3). (Recall that $\dim\Om^3_k=\binom{k-1}{2}$.)
\msn
{\bf 1.8.~Remark.} It is known that $\al_Z>\tfrac{2}{d}$ for any reduced projective hypersurface $Z\subset\PP^{n-1}$ of degree $d$, assuming that $f$ is not a polynomial of $n'$ variables with $n'<n$ and $n\ges 3$. This follows from \cite[Corollary~3.6]{dFEM}. (Note that $\al_Z$ coincides with the {\it log canonical threshold} of $Z\subset\PP^{n-1}$ as is well known.) This assertion does not seem to follow easily from the semicontinuity argument. Indeed, if one considers a one-parameter family $tf+sg$ where $[t\1{:}\1s]\in\PP^1$ and $g$ is an appropriate homogeneous polynomial of 2 variables with degree $d$, then it is entirely unclear why the Milnor number of a given singular point of the curve can be made {\it stable} for $|s/t|$ sufficiently small by choosing $g$ appropriately. This is quite different from the case of adding higher monomials using the finite determinacy property of isolated hypersurface singularities.
\msn
{\bf 1.9.~Remark.} In the line arrangement case with $n=3$, the following two assertions are equivalent to each other:
$$\mu_k=\tau_Z\,\,\,\,\,\h{if}\,\,\,\,\,k\ges 2d-1,
\leqno(1.9.1)$$
$$\mu^{(\infty)}_k=0\,\,\,\,\,\h{if}\,\,\,\,\,k\ges 2d-1.
\leqno(1.9.2)$$
Indeed, we have by \cite[Theorem~3 and Corollary~3]{DiSa2}
$$\aligned&{}'\!\ddd^{(1)}:N_{q+d}\to M_q\,\,\,\h{is {\it injective,} and consequently}\\ &{}'\!\ddd^{(r)}:N^{(r)}_{q+d}\to M^{(r)}_{q+d-rd}\,\,\,\h{\it vanishes}\,\,(\forall\,r\ges2),\,\,\,\h{if}\,\,\,\,q\ges 2d-1,\endaligned
\leqno(1.9.3)$$
$$\nu_{k+d}=\tau_Z\,\,\,\,\,\h{if}\,\,\,\,\,k\ges 2d-1.
\leqno(1.9.4)$$
Note that (1.9.2) is equivalent to \cite[Theorem~1]{bha} in the line arrangement case by \cite[Theorem~2]{bfun}. For (1.9.1), see \cite{DIM}.
\msn
{\bf 1.10.~Remark.} We have the {\it inequality} $F^{\ssb}\ne P^{\ssb}$ on $H^{n-1}(\Ff,\C)_{\lambda}$ of any indecomposable essential reduced central hyperplane arrangements in $\C^n$ ($n\ges 3$) at least for $\lambda=\exp(2\pi i/d)$ or $\exp(-2\pi i/d)$ (or more generally, if $\lambda$ is not an eigenvalues of the Milnor monodromy of $f$ at any $x\in\C^n\setminus\{0\}$). Indeed, the vanishing cycle sheaf $\varphi_{f,\lambda}\C[n-1]$ is supported at the origin for these $\lambda$, since $f$ is essential. So it can be identified with the Milnor cohomology $H^{n-1}(\Ff,\C)_{\lambda}$, and moreover the direct sum $\varphi_{f,\lambda}\C[n-1]\oplus\varphi_{f,\overline{\lambda}}\C[n-1]$ is identified with a real Hodge structure of weight $2$ (since $N=0$) so that we can get a symmetry of its Steenbrink spectral numbers with center $\frac{n}{2}$ as in the isolated singularity case. Assume $F^{\ssb}=P^{\ssb}$ on this direct sum. By the indecomposability assumption we have the nonvanishing of the Euler characteristic $\chi(U)$ (see \cite{STV}), hence $1-\frac{1}{d}$ must be a spectral number, since it must be {\it strictly smaller} than $2-\frac{1}{d}$ by \cite[Theorem~1]{bha}. Then $n-1+\frac{1}{d}$ must be also a spectral number by the symmetry, although it is not allowed by \cite[Theorem~1]{bha} as is explained above (since $n\ges 3$). So we get $F^{\ssb}\ne P^{\ssb}$ on the direct sum.
\bs\bs
\vbox{\centerline{\bf 2. Proof of main theorems}
\bsn
In this section we prove Theorems~2 and 3 together with Propositions~1.} 
\msn
{\bf 2.1.~Roots of Bernstein-Sato polynomials.} In the notation of the introduction, set
$$\al_f:=\min\R_f\,,\q\al_Z:=\min\R_Z\,,\q\alt_f:=\min\Rt_f\,,\q\alt_Z:=\min\Rt_Z\,,$$
where $\Rt_f\subset\R_f$ is the set of roots of the {\it microlocal\,} (that is, reduced) Bernstein-Sato polynomial $\bt_f(s)=b_f(s)/(s+1)$ up to a sign (see \cite{mic}), and similarly for $\Rt_Z=\mcup_{z\in Z}\,\Rt_{h_z,z}$ in the notation of the introduction. Set
$$\alt'_Z:=\min\Rt_{Z'},$$
where $Z'$ is a general hyperplane section of $Z$. We have
$$\al_f=\min\bl(\al_Z,\nod\br)\les\alt_f\les\alt_Z\les\alt'_Z.
\leqno(2.1.1)$$
The first equality is a consequence of \cite[Theorem 2.2]{bfun} (see also\cite[Theorem 4.8]{wh}) together with the independence of Bernstein-Sato polynomials on the choice of defining functions as is noted at the beginning. The first inequality is trivial, since $\al_f=\min(\alt_f,1)$. The second one follows from the independence of a choice of a defining function mentioned just above. The last one can be shown by using \cite[Lemma~4.2]{DMST}. It is not necessarily easy to determine $\alt_f$ in the case $\alt_f>1$, although $\al_f$ is easily obtained by using $\al_Z$ and $\nod$ via the first equality of (2.1.1).
\sk
We have by \cite[Theorem~0.4]{mic}
$$\max\Rt_f\les n-\alt_f,\q\max\Rt_Z\les n-1-\alt_Z,\q\max\Rt_{Z'}\les n-2-\alt'_Z.
\leqno(2.1.2)$$
This implies that $\alt_f\les\tfrac{n}{2}$, $\alt_Z\les\tfrac{n-1}{2}$, $\alt'_Z\les\tfrac{n-2}{2}$. Hence (2.1.2) holds with $\Rt$ replaced by $\R$.
\msn
{\bf Remark.}
Assuming $\dim{\rm Sing}\,Z\les 1$, we have
$$\aligned\mu^{(\infty)}_k&=0\q\h{if}\,\,\,\,\kod>n-\alt_f,\\ \nu^{(\infty)}_k&=0\q\h{if}\,\,\,\,\kod>n-\alt_Z,\\ \rho^{(\infty)}_k&=0\q\h{if}\,\,\,\,\kod>n-\alt'_Z.\endaligned
\leqno(2.1.3)$$
The first and last assertions follow from (2.1.1--2) together with Theorem~1 and Proposition~1. However the proof of the middle is quite complicated, and its proof is omitted, since it is not needed in this paper.
\msn
{\bf 2.2.~Pole order spectrum.} For a homogeneous polynomial $f$ of $n$ variables and degree $d$, set $F_{\!f}:=f^{-1}(1)$, the Milnor fiber of $f$. Define
$$\aligned&\q\q\q\q\Sp_P^j(f):=\msum_{\al\in\Q}\,{}^P\!n^j_{f,\al}\,t^{\al}\q\q(j\in\N)\\
\h{with}\q\q&{}^P\!n^j_{f,\al}:=\dim\Gr^p_P\Ht^{n-1-j}(F_{\!f},\C)_{\ee(-\al)}\q\bl(\1 p=[n-\al]\br).\endaligned
\leqno(2.2.1)$$
Here $\ee(-\al)$ is as in Theorem~1, and $P$ is the {\it pole order filtration} associated with the pole order spectral sequence, see \cite{DiSa2}. We have by definition
$$\mu^{(\infty)}_k={}^P\!n^0_{f,k/d},\q\nu^{(\infty)}_k={}^P\!n^1_{f,k/d},\q\rho^{(\infty)}_k={}^P\!n^2_{f,k/d}\q\q(k\in\N).
\leqno(2.2.2)$$
\msn
{\bf 2.3.~Proof of Proposition~1.} Let $\Om^{\prime\ssb}$ be the complex of global sections of algebraic differential forms on a general hyperplane of $\C^n$ containing the origin. There is a natural restriction morphism
$$\iota:\Om^{\ssb}\to\Om^{\prime\ssb},
\leqno(2.3.1)$$
which is compatible with $\ddd$ and $\wedge$. So it induces a morphism between the microlocal pole order spectral sequences for $f$ and $g$, where $g$ is the restriction of $f$ to a hyperplane so that $Z'=\{g=0\}$ in the notation of the introduction.
\sk
The assertion (10) then follows from a similar injectivity assertion in \cite[Theorem~5.3]{DiSa2} for hyperplane sections together with a kind of ``torsion-freeness`` of $Q$ (after taking the direct image by a projection to $\C^2$), more precisely, the equality $Q=Q'''$ in (3.4.1) below, where the local cohomology functor $H^0_I$ in (3.1.3) is used. By this ``torsion-freeness", it is enough to apply the above restriction argument to {\it sufficiently general\,} members of a general $1$-parameter family of hyperplanes of $\PP^{n-1}$ (such that each irreducible component of ${\rm Sing}\,Z$ is not contained in general members of it as is noted after Theorem~2). By condition~(GH), $\Rt_{Z'}$ coincides with the {\it spectral numbers} of $Z'$ (see for instance \cite[Section 1.9]{wh}), and we have by (2.1.2)
$$\max\Rt_{Z'}\les n-2-\alt'_Z.$$
There is a shift of indices $k$ by $d$ coming from \cite[Theorem~5.3]{DiSa2}, and we have another shift by the restriction morphism $\iota$ in (2.3.1), since the shift by the ambient dimension $n$ is used for the correspondence between the pole order filtration and the pole order spectrum in (2.2.1), see also (1) in Theorem~1. This finishes the proof of Proposition~1.
\msn
{\bf 2.4.~Proof of Theorem~2.} It is well known that there are local systems $\LL_{\la}$ of rank 1 on $U:=\PP^{n-1}\setminus Z$ calculating $H^{\ssb}(F_{\!f},\C)_{\la}$, that is,
$$H^j(F_{\!f},\C)_{\la}=H^j(U,\LL_{\la})\q\q(\forall\,\la\in\mmu_d),
\leqno(2.4.1)$$
see for instance \cite[1.4.2]{BuSa}, \cite{wh}, etc. This implies that
$$\msum_i\,(-1)^iH^i(\Ff,\C)_{\la}=\chi(U).$$
\sk
Using the first and last assertions in (2.1.3), we then get
$$\chi^{(\infty)}_{f,\la}(\beta)\ges(-1)^{n-1}\chi(U)\q(\forall\,\la\in\mmu_d).
\leqno(2.4.2)$$
since $\beta\les\lceil\alt'_Z\rceil^{/d}$ by condition~(5) and (2.1.1). Proposition~1 together with the last inequality for $\beta$ implies that
$${}''\!\ddd^{(m)}:Q^{(m)}_k\to N^{(m)}_{k-md}\,\,\,\,\h{vanishes if}\,\,\,\,\kod>n-\beta,\,\,m\ges 2.
\leqno(2.4.3)$$
We then get inequalities
$$\chi^{(m+1)}_{f,\la}(\beta)\les\chi^{(m)}_{f,\la}(\beta)\q\q(\forall\,\la\in\mmu_d,\,\,m\ges 2),
\leqno(2.4.4)$$
where {\it strict inequalities\1} occur if and only if there are some nonzero differentials
$${}'\!\ddd^{(m)}:N^{(m)}_k\to M^{(m)}_{k-md}\q\h{with}\q\kod\in(n-\beta,\,n+m-\beta].
\leqno(2.4.5)$$
Note that there may be differentials
$$\aligned{}'\!\ddd^{(r)}&:N^{(m)}_k\to M^{(m)}_{k-md}\q\h{with}\q \kod\les n-\beta,\\{}''\!\ddd^{(r)}&:Q^{(m)}_k\to N^{(m)}_{k-md}\q\h{with}\q \kod\les n-\beta,\endaligned
\leqno(2.4.6)$$
{\it which do not contribute to the change from $\chi^{(m)}_{f,\la}(\beta)$ to $\chi^{(m+1)}_{f,\la}(\beta)$ even if they are nonzero}. We should keep them in mind for the proof of our {\it partial $E_r$-degeneration} of the spectral sequence as in (8). We can forget about the differentials
$${}'\!\ddd^{(m)}:N^{(m)}_k\to M^{(m)}_{k-rd}\q\h{with}\q \kod>n+m-\beta,$$
since these are irrelevant to the partial degeneration, and do not contribute to the change from $\chi^{(m)}_{f,\la}(\beta)$ to $\chi^{(m+1)}_{f,\la}(\beta)$.
\sk
By condition~(6) together with (2.4.2) and (2.4.4), we have the following inequalities for any $\la\in\mmu_d$:
$$(-1)^{n-1}\chi(U)\les\chi^{(\infty)}_{f,\la}(\beta)\les\chi^{(r)}_{f,\la}(\beta)\les(-1)^{n-1}\chi(U),
\leqno(2.4.7)$$
since $r\ges 2$. All the inequalities then become equalities, and we get in particular
$$\chi^{(\infty)}_{f,\la}(\beta)=\chi^{(r)}_{f,\la}(\beta),
\leqno(2.4.8)$$
together with the vanishing of (2.4.5) for $m\ges r$. The differential ${}'\!\ddd^{(m)}$ in (2.4.6) vanishes for $m\ges r$ by condition~(7), and similarly for ${}''\!\ddd^{(m)}$ about the last assertion of Theorem~2. This finishes the proof of Theorem~2.
\msn
{\bf 2.5.~Remarks.} (i) Condition~(AT) in the introduction means that, at any point $z$ of a sufficiently small Zariski open subset of $Z$, there is a germ of an analytically defined smooth morphism $\pi_z:(\PP^{n-1},z)\to(\C^{n-2},0)$ such that $Z$ is locally the pull-back of a hypersurface of $(\C^{n-2},0)$ by $\pi_z$. Condition~(GH) seems to imply (AT) except for certain special cases (for instance, where the general hyperplane section of $Z$ has a weighted homogeneous isolated singularity such that any $\mu$-constant deformations are still weighted homogeneous, that is, the base space of the versal deformation \cite{KaSc}, \cite{Tj} has only non-negative weight part, see for instance Example~5.9 below).
\ms
(ii) Condition~(GH) in Theorem~2 may be replaced by a slightly weaker hypothesis that the condition is satisfied only for general members of a certain 1-parameter family of hyperplane sections such that each irreducible component of ${\rm Sing}\,Z$ is not contained in general members of it, see (2.3).
\ms
(iii) Condition~(GH) holds trivially when $n=3$ (here $Z$ may be a {\it non-reduced\1} curve). For $n=4$, is satisfied by reduced hyperplane arrangements and also by locally positively weighted homogeneous reduced divisors, see \cite[Proposition 2.4]{CNM}.
\msn
{\bf 2.6.~Proof of Theorem~3.} The argument is essentially the same as the proof of Theorem~2. We have the inequalities
$$-\chi(U)\les\cho^{(\infty)}_{f,\la}(\beta)\les\cho^{(r)}_{f,\la}(\beta)\les-\chi(U),
\leqno(2.6.1)$$
where the first inequality follows from the assertion~(15) and Proposition~1, the second one is shown by an argument similar to the proof of Theorem~2, and the last one is condition~(12) in Theorem~3. All the inequalities then become {\it equalities}.
\sk
As a consequence, condition~(14) must hold with $\nu^{(r)}_k$ replaced by $\nu^{(\infty)}_k$, since the first inequality of (2.6.1) becomes a strict inequality otherwise. The assertion for the second inequality together with Proposition~1 then implies that condition~(14) hold for $\nu^{(r)}_k$, since the differential ${}'\!\ddd^{(m)}$ decreases the degree by $-md$ (and $m\ges r\ges 2$).
\sk
By a similar argument, we get the partial $E_r$-degeneration of the pole order spectral sequence for the highest cohomology with $\kod\les 2-\twd\,$, using Proposition~1 together with condition~(13). The argument is similar for the assertion corresponding to the last assertion of Theorem~2. This finishes the proof of Theorem~3.
\msn
{\bf Remark.} It is unclear whether we have the equalities
$$\cho^{(2)}_{f,\la}\bl(\twd\br)=\chi^{(2)}_{f,\la}\bl(\twd\br)\q\q\bl(\la\in\mmu_d\br),$$
since the vanishing of $\mu^{(2)}_k$ for $k>m$ is nontrivial, and conditions (14--15) imply only that
$$\cho^{(\infty)}_{f,\la}\bl(\twd\br)=\chi^{(\infty)}_{f,\la}\bl(\twd\br)\q\q\bl(\la\in\mmu_d\br).$$
Here we cannot control the differentials ${}'\!\ddd^{(m)}:N^{(m)}_{k+md}\to M^{(m)}_k$ for $\kod\in[2-\ond,4-\twd]$, $m\ges 2$.
\bs\bs
\vbox{\centerline{\bf 3. Self-duality isomorphisms}
\bsn
In this section we prove a structure theorem together with self-duality isomorphisms for the Koszul cohomologies $M,N,Q$. These imply a kind of ``torsion-freeness" of $Q$.}
\msn
{\bf 3.1.~Spectral sequences.} Let $R=\C[x_1,\dots,x_n]$, and $\Om^j=\Gamma(\C^n,\Om_{\C^n}^j)$ ($j\in\Z$) as in (1.1). For a bounded complex of finitely generated graded $R$-modules $M^{\ssb}$, set
$$\aligned\D(M^{\ssb})&:=\RR{\rm Hom}_R(M^{\ssb},\Om^n[n]),\\ D_i(M^{\ssb})&:=H^{-i}\bl(\D(M^{\ssb})\br)={\rm Ext}_R^{n-i}(M^{\ssb},\Om^n)\q(i\in\Z).\endaligned$$
We have the following spectral sequences in the abelian category of graded $R$-modules
$$E_2^{p,q}=D_{-p}(H^{-q}M^{\ssb})\Longrightarrow D_{-p-q}(M^{\ssb}),
\leqno(3.1.1)$$
as is well-known, see for instance \cite{De} (and also \cite{DiSa2}).
\sk
We apply these to the shifted Koszul complex ${}^s\!K_f^{\ssb}:=(\Om^{\ssb},\df\wedge)[n]$ in \cite{DiSa2}. Here the grading of ${}^s\!K_f^{\ssb}$ is shifted so that the differential $\df\wedge$ preserves it, that is,
$${}^s\!K_{f,k}^j=\Om^{n+j}_{k+jd}\q(j,k\in\Z),$$
and we have the self-duality
$$\D({}^s\!K_f^{\ssb})={}^s\!K_f^{\ssb}(nd).$$
So (3.1.1) gives the spectral sequence
$$E_2^{p,q}=D_{-p}\bl(H^{-q}({}^s\!K_f^{\ssb})\br)\Longrightarrow H^{p+q}({}^s\!K_f^{\ssb})(nd),
\leqno(3.1.2)$$
with
$$H^0({}^s\!K_f^{\ssb})=M,\q H^{-1}({}^s\!K_f^{\ssb})=N,\q H^{-2}({}^s\!K_f^{\ssb})=Q,$$
and $H^i({}^s\!K_f^{\ssb})=0$ otherwise (since we assume $\dim{\rm Sing}\,f^{-1}(0)\les 2$).
\sk
Let $\mm:=(x_1,\dots,x_n)\subset R$, the graded maximal ideal of $R$. Let $I$ be the reduced graded ideal of $R$ corresponding to a sufficiently large finite subset $\Sigma\subset{\rm Sing}\,Z$ containing the singular locus of $({\rm Sing}\,Z)_{\rm red}$. (Note that $\dim{\rm Sing}\,Z\les 1$ by assumption.) Set
$$M':=H_{\mm}^0(M),\q M'':=H_I^0(M)/H_{\mm}^0(M),\q M''':=M/H_I^0(M),
\leqno(3.1.3)$$
(similarly for $N'$, $Q'$, etc.)
These are independent of $I$ as long as $I$ is sufficiently large. Indeed, let $G$ be the filtration on $M$ {\it by codimension of support,} that is,
$$G^pM=\bl\{\om\in M\mid{\rm codim}_{S_f}\,{\rm Supp}\,R\1\om\ges p\br\}\q\h{with}\q S_f:={\rm Sing}\,f^{-1}(0).$$
Then
$$G^0M=M,\q G^1M:=H_I^0(M),\q G^2M:=H_{\mm}^0(M),\q G^3M=0,$$
and
$$\Gr_G^0M=M''',\q\Gr_G^1M=M'',\q\Gr_G^2M=M'.$$
This filtration $G$ is trivial (that is, $G^1M=G^2M=0$) if $Z$ is strongly free (see Remark~(iii) after (4.1) below). Its converse is unclear, since $M'''_{\rm def}$ may be nonzero.
\sk
We have the spectral sequence associated with the filtration $G$:
$$E_1^{p,q}=D_{-p-q}(\Gr_G^{-p}M)\Longrightarrow D_{-p-q}(M),
\leqno(3.1.4)$$
(similar assertions hold with $M$ replaced by $N,Q$), see also \cite{BrHe}, \cite{Ei1}, etc.
\sk
Let $\pi_i$ be a generic projection form $\C^n$ to $\C\1^i$ ($i=1,2$) inducing {\it finite} morphisms
$$\overline{\pi}_1:\Sigma\to\PP^0,\q\overline{\pi}_2:{\rm Sing}\,Z\to\PP^1.
\leqno(3.1.5)$$
Let $(M'')^{\sim}$, $(M''')^{\sim}$ be the coherent sheaves on $\C^n$ corresponding to the $R$-modules $M''$, $M'''$ as in \cite{Gr2} or \cite[II, Corollary 5.5]{Ha}. Then
$$\h{$(\pi_1)_*(M'')^{\sim}$ and $(\pi_2)_*(M''')^{\sim}|_{\C^2\setminus\{0\}}\,$ are free sheaves of finite rank.}
\leqno(3.1.6)$$
For the assertion about $(M''')^{\sim}$, note that any torsion-free coherent sheaf on $\PP^1$ is locally free, and is isomorphic to a direct sum of line bundles $\OO_{\PP^1}(a_i)$ with $a_i\in\Z$, see \cite{Gr1}.
Here we use a well-known relation between graded modules over an affine ring and coherent sheaves on the associated projective space, see \cite{Gr2} or \cite[II, Exercise 5.9(c)]{Ha}.
\sk
Let $M'''_{\rm max}$ be the graded $R$-module such that
$$(M'''_{\rm max})^{\sim}=(j_0)_*j_0^*(M''')^{\sim},$$
where $j_0:\C^n\setminus\{0\}\into\C^n$ is the natural inclusion. Then
$$\h{$(\pi_2)_*(M'''_{\rm max})^{\sim}\,$ is a free sheaf on $\,\C^2$.}
\leqno(3.1.7)$$
We have a short exact sequence of graded $R$-modules
$$0\to M'''\to M'''_{\rm max}\to M'''_{\rm def}\to 0,
\leqno(3.1.8)$$
where $M'''_{\rm def}$ is finite dimensional over $\C$, and is defined by the exact sequence.
\sk
Let $(M''')^{\sim}_{\PP^{n-1}}$ be the coherent sheaf on $\PP^{n-1}$ associated with the graded $R$-module $M'''$ as in {\it loc.\,cit.} Then there are canonical isomorphisms
$$M'''_{\rm max}=\Gamma_*(M''')^{\sim}_{\PP^{n-1}}:=\mopl_{i\in\Z}\,H^0\bl(\PP^{n-1},(M''')^{\sim}_{\PP^{n-1}}(i)\br),
\leqno(3.1.9)$$
$$M'''_{\rm def}=H^1_{\mm}(M'''),
\leqno(3.1.10)$$
where $\mm\subset R$ is the graded maximal ideal as in (3.1.3), and the exact sequence (3.1.8) can be identified with the long exact sequence associated with local cohomology and open direct image, see for instance \cite[Proposition 2.1.5]{Gr3}, \cite[Corollary A1.12]{Ei2}. (Similar assertions hold with $M$ replaced by $N$, $Q$.)
\sk
We have the following.
\msn
{\bf 3.2.~Lemma.} {\it In the above notation, we have
$$\aligned\,\,&D_i(M')=0\,\,\,(i\ne 0)\,\,\,\bl(\h{that is,}\,\,\,\D(M')=D_0(M')\br),\,\,\,D_0^2(M')=M',\\&(\h{and a similar assertion holds with $\,M'\,$ replaced by $\,M'''_{\rm def})$}, \endaligned
\leqno(3.2.1)$$
$$D_i(M'')=0\,\,\,(i\ne 1)\,\,\,\bl(\h{that is,}\,\,\,\D(M'')=D_1(M'')[1]\br),\,\,\,D_1^2M''=M'',
\leqno(3.2.2)$$
$$\aligned&D_i(M''')=0\,\,\,(i\ne 1,2),\q D_i(M'''_{\rm max})=0\,\,\,(i\ne 2),\\&D_1(M''')=D_0(M'''_{\rm def}),\,\,\,\,D_2(M''')=D_2(M'''_{\rm max}),\\&D_0D_1(M''')=M'''_{\rm def},\,\,\,D_2^2(M''')=D_2^2(M'''_{\rm max})=M'''_{\rm max}.
\endaligned\leqno(3.2.3)$$
Moreover similar assertions hold with $M$ replaced by $N$, $Q$.}
\msn
{\it Proof.} This follows from (3.1.6--8) together with the compatibility of the duality with the direct image by $\pi_1$ or $\pi_2$ (where the support of the sheaf is finite over $\C^1$ or $\C^2$). For (3.2.3) we use the following long exact sequence induced from (3.1.8):
$$\to D_i(M'''_{\rm def})\to D_i(M'''_{\rm max})\to D_i(M''')\to D_{i-1}(M'''_{\rm def})\to D_{i-1}(M'''_{\rm max})\to$$
This finishes the proof of Lemma~(3.2).
\ms
Combined with the spectral sequence (3.1.4), this implies the following.
\msn
{\bf 3.3.~Proposition.} {\it In the notation of $(3.1)$, we have the isomorphisms
$$D_0(M)=D_0(M'),
\leqno(3.3.1)$$
$$D_2(M)=D_2(M''')=D_2(M'''_{\rm max}),
\leqno(3.3.2)$$
together with the short exact sequence
$$0\to D_0(M'''_{\rm def})\to D_1(M)\to D_1(M'')\to 0.
\leqno(3.3.3)$$
Moreover similar assertions hold with $M$ replaced by $N$, $Q$.}
\msn
{\it Proof.} The $E_1$-terms $E_1^{p,q}$ of the spectral sequence (3.1.4) are given by
$$\begin{array}{cccccccc}D_0(M'),\\ D_1(M'),&D_0(M''),&\\ D_2(M'),&D_1(M''),&D_0(M''')\\ &D_2(M''),&D_1(M''')\\ &&D_2(M''')\end{array}$$
where the other terms vanish, since
$$D_i(M')=D_i(M'')=D_i(M''')=0\q(i>2),$$
see Lemma~(3.2). Using further the other vanishing assertions in the lemma, we get the $E_1$-degeneration of the spectral sequence (3.1.4), and the assertion follows. This finishes the proof of Proposition~(3.3).
\ms
Using the spectral sequence (3.1.2), these imply the following.
\msn\vbox{\nin
{\bf 3.4.~Theorem.} {\it In the notation of $(3.1)$, we have the isomorphisms of graded $R$-modules
$$D_2(M)=Q(nd)=Q'''_{\rm max}(nd),\q\q Q'=Q''=Q'''_{\rm def}=0,
\leqno(3.4.1)$$
$$D_0(M'''_{\rm def})=D_1(M''')=N'(nd),
\leqno(3.4.2)$$
$$D_1(M'')=N''(nd),
\leqno(3.4.3)$$
$$D_2(M'''_{\rm max})=Q'''_{\rm max}(nd),\q\q D_2(N'''_{\rm max})=N'''_{\rm max}(nd).
\leqno(3.4.4)$$
together with the exact sequence}
$$0\to N'''_{\rm def}(nd)\to D_0(M')\buildrel{\phi}\over\longrightarrow M'(nd)\to D_0(N'''_{\rm def})\to 0.
\leqno(3.4.5)$$}
\msn
{\it Proof.} Using the isomorphisms after (3.1.2) together with Proposition~(3.3), the $E_2$-terms $E_2^{p,q}$ of the spectral sequence (3.1.2) are given by
$$\begin{array}{cccccccc}D_2(Q'''_{\rm max}),&D_1(Q),&D_0(Q')\\ D_2(N'''_{\rm max}),&D_1(N),&D_0(N')\\ D_2(M'''_{\rm max}),&D_1(M),&D_0(M')\end{array}$$
where the other terms vanish. The differential $\ddd_2$ vanishes except for
$$\ddd_2:D_2(Q'''_{\rm max})\to D_0(N'),\q\ddd_2:D_2(N'''_{\rm max})\to D_0(M').$$
So the spectral sequence (3.1.2) degenerates at $E_3$, and we get the desired isomorphisms and exact sequences using Lemma~(3.2) and Proposition~(3.3). Related to the differential $\ddd_2$ of the spectral sequence (3.1.2), for instance, we have the short exact sequences
$$0\to M'''(nd)\to D_2(Q'''_{\rm max})\buildrel{\ddd_2\,\,}\over\to D_0(N')\to 0,
\leqno(3.4.6)$$
$$0\to{\rm Coker}\bl(D_2(N'''_{\rm max})\buildrel{\ddd_2\,\,}\over\to D_0(M')\br)\to M'(nd)\to D_0(N'''_{\rm def})\to 0,
\leqno(3.4.7)$$
together with the isomorphism
$$N'''(nd)={\rm Ker}\bl(D_2(N'''_{\rm max})\buildrel{\ddd_2\,\,}\over\to D_0(M')\br).
\leqno(3.4.8)$$
Here (3.4.6) implies the first isomorphism of (3.4.4) together with (3.4.2).
As for (3.4.7--8), these are equivalent to the exact sequence
$$0\to N'''(nd)\to D_2(N'''_{\rm max})\buildrel{\ddd_2\,\,}\over\to D_0(M')\to M'(nd)\to D_0(N'''_{\rm def})\to 0,
\leqno(3.4.9)$$
which is further equivalent to the exact sequence (3.4.5) together with
$$0\to N'''(nd)\to D_2(N'''_{\rm max})\to N'''_{\rm def}(nd)\to 0,
\leqno(3.4.10)$$
where the latter implies the last isomorphism of (3.4.4).
The argument is easier for other isomorphisms. This finishes the proof of Theorem~(3.4).
\msn
{\bf Remark.} The exact sequence (3.4.5) shows that the self-duality of $M'$ does not necessarily hold unless $N'''_{\rm def}=0$.
\msn
{\bf 3.5.~Numerical formulas.} For positive integers $m$, we have the sequences $p^{(m)}=\{p^{(m)}_k\}_{k\in\Z}$ defined by
$$p^{(m)}_k:=\begin{cases}\tbinom{k+m-1}{m-1}&\h{if}\,\,\,\,k\ges 0,\\ \q\,\,\,0&\h{if}\,\,\,\,k<0.\end{cases}$$
These are the dimensions of the vector spaces of homogeneous polynomials of $m$ variables with degree $k$. For $m=1,2,3$ and $k\ges 0$, we have
$$p^{(1)}_k=1,\q\q\q p^{(2)}_k=k+1,\q\q\q p^{(3)}_k=\tfrac{(k+1)(k+2)}{2}.
\leqno(3.5.1)$$
\sk
For a sequence $q=\{q_k\}_{k\in\Z}$, define
$${\rm Diff}(q)_k:=q_k-q_{k-1}\q\q(k\in\Z).$$
Then we have for $m\ges 1$
$${\rm Diff}(p^{(m)})_k=p^{(m-1)}_k\q\q(k\in\Z),
\leqno(3.5.2)$$
$${\rm Diff}^m(p^{(m)})_k=p^{(0)}_k:=\begin{cases}1\q&\h{if}\,\,\,\,k=0,\\0&\h{if}\,\,\,\,k\ne 0.\end{cases}
\leqno(3.5.3)$$
In this paper we need also this $p^{(0)}_k$, since we do not consider the graded $R$-modules modulo lower degrees as in \cite[I, Proposition 7.3]{Ha}.
\sk
Set
$$\mu'_k:=\dim M'_k,\q\mu''_k:=\dim M''_k,\q\mu'''_k:=\dim M'''_k,$$
$$\mu'''_{{\rm max},\,k}:=\dim M'''_{{\rm max},\,k},\q\mu'''_{{\rm def},\,k}:=\dim M'''_{{\rm def},\,k},$$
and similarly for $\nu'''_{{\rm def},\,k}$, $\rho'''_{{\rm def},\,k}$, etc.
By (3.1.6--7) there are integers
$$a''_i,\,a'''_j\in\Z\,\,\,\,\bl(i\in[1,r''],\,\,j\in[1,r''']\br),$$
with $r''$, $r'''$ non-negative integers and such that
$$\mu''_k=\msum_{i=1}^{r''}\,p^{(1)}_{k-a''_i},\q\q\mu'''_{{\rm max},\,k}=\msum_{j=1}^{r'''}\,p^{(2)}_{k-a'''_j}\q\q(k\in\Z),
\leqno(3.5.4)$$
where $r''$, $r'''$ are respectively called the {\it rank\1} of $M''$, $M'''$, and similarly for $\nu''_k$, $\nu'''_{{\rm max},\,k}$, and also for $\rho_k=\rho'''_{{\rm max},\,k}$, with $a''_i,a'''_j,r'''$ replaced by $b''_i,b'''_j,2\1r'''$, and $a'''_j$ replaced by $c'''_j$. Here
$$\aligned r''&={\rm rank}\,M''={\rm rank}\,N'',\\ r'''&={\rm rank}\,M'''={\rm rank}\,Q'''=\tfrac{1}{2}\1{\rm rank}\,N''',\endaligned
\leqno(3.5.5)$$
by Theorem~(3.4) and (4.2.5) below. (We can show that $r'''$ is at most the Tjurina number $\tau_{Z'}$ of a general hyperplane section $Z'$ of $Z$. They coincide if condition~(AT) in the introduction is satisfied.) The $\mu'_k$, $\mu'''_{{\rm def},\,k}$, etc.\ are finite linear combinations of the $p^{(0)}_{k-i}$ ($i\in\Z$) with positive coefficients.
\sk
Theorem~(3.4) implies the following.
\msn
{\bf Corollary~3.6.} {\it In the above notation, we have the following equalities for any $k\in\Z$:}
$${\rm Diff}^2(\mu'''_{\rm max})_k={\rm Diff}^2(\rho)_{nd+2-k}\,,
\leqno(3.6.1)$$
$${\rm Diff}^2(\nu'''_{\rm max})_k={\rm Diff}^2(\nu'''_{\rm max})_{nd+2-k}\,,
\leqno(3.6.2)$$
$${\rm Diff}(\mu'')_k={\rm Diff}(\nu'')_{nd+1-k}\,,
\leqno(3.6.3)$$
$$\mu'''_{{\rm def},\,k}=\nu'_{nd-k}\,,
\leqno(3.6.4)$$
$$\mu'_k-\nu'''_{{\rm def},\,nd-k}=\mu'_{nd-k}-\nu'''_{{\rm def},\,k}\,\ges \,0,
\leqno(3.6.5)$$
$$\rho_k=\rho'''_k=\rho'''_{{\rm max},\,k},\q\rho'_k=\rho''_k=\rho'''_{{\rm def},\,k}=0\,.
\leqno(3.6.6)$$
\msn
{\it Proof.} These follow from Theorem~(3.4). For instance, we can deduce (3.6.5) from the exact sequence (3.4.5). For (3.6.1--3) we also use the compatibility of the duality with the direct images by $\pi_2$, $\pi_1$ in (3.1) together with (3.5.3--4). Here we have the shift of index by $2$ or $1$ on the right-hand side of (3.6.1--3), since the grading is shifted by $m$ on the dualizing complex on $\C^m$ ($m=1,2$) as is seen in (3.1) for the case $m=n$. This finishes the proof of Corollary~(3.6).
\msn
{\bf 3.7.~Remarks.} (i) The assertion (3.6.3) is equivalent to the following:
$$\h{$\mu''_k+\nu''_{nd-k}\,$ is constant ($k\in\Z$).}
\leqno(3.7.1)$$
\sk
(ii) We can determine the $\mu'_k$, $\mu''_k$ by calculating the {\it saturation} of the ideal $(\dd f)\subset R$ for the maximal ideal and for an {\it certain\1} ideal $E\subset R$ using a computer program like Macaulay2 or Singular, see \cite{GrSt}, \cite{DGPS}. (When $n=4$, $E$ is the {\it product\1} of the {\it associated primes} of the Jacobian ideal $(\dd f)\subset R$ with {\it codimension} 3, which can be {\it embedded primes}.) In the case of Example~(5.8) below and Macaulay2, one may type M2 and press RETURN on a terminal of Unix (or Mac, etc.) with Macaulay2 installed, and then copy and past the following (after some corrections if necessary, see Note below):
\ms
\vbox{\small\sf\verb#R=QQ[x,y,z,w]; f=x*y*z*w*(x+y+z)*(y-z+w);#
\sk
\verb#J=ideal(jacobian ideal(f)); K=saturate(J,ideal(x,y,z,w));#
\sk
\verb#E=ideal(y,z,x*w); L=saturate(J,E); d=first degree f;#
\sk
\verb#mxd=4*d; A=QQ[v]/(v^mxd);#
\sk
\verb#mupA=v^4*sub(hilbertSeries(K/J,Order=>mxd),vars A)#
\sk
\verb#musA=v^4*sub(hilbertSeries(L/K,Order=>mxd),vars A)#
\sk
\verb#mutA=v^4*sub(hilbertSeries(R/L,Order=>mxd),vars A)#
\sk
\verb#muA=mupA+musA+mutA#}
\skn
In the case of Singular, one may do, for instance, as follows (where one gets a different kind of output that may be easier to see):
\ms
\vbox{\small\sf\verb#LIB "elim.lib";#
\sk
\verb#ring R = 0, (x,y,z,w), dp; int i, d; intvec k, mu, mup, mus, mut;#
\sk
\verb#poly f=x*y*z*w*(x+y+z)*(y-z+w);#
\sk
\verb#d=deg(f); ideal J=jacob(f); ideal I=sat(J,maxideal(1))[1];#
\sk
\verb#ideal E=(y,z,x*w); ideal K=sat(J,E)[1];#
\sk
\verb#for (i=1; i<4*d-3; i++) {#
\sk\q
\verb#mu[i]=size(kbase(std(J),i-1)); mup[i]=mu[i]-size(kbase(std(I),i-1));#
\sk\q
\verb#mut[i]=size(kbase(std(K),i-1)); mus[i]=mu[i]-mup[i]-mut[i]; k[i]=i+3; }#
\sk
\verb#sprintf(" k : %s",k); sprintf("mup: %s",mup); sprintf("mus: %s",mus);#
\sk
\verb#sprintf("mut: %s",mut); sprintf("mu : %s",mu);#}
\msn
{\bf Note.} If one copies these from a pdf file, RETURN {\it must be added if it is erased.} Sometimes $\{$ and $\}$ may be changed to f and g, depending on the viewer, and these must be corrected.
\msn
{\bf Caution 1.} The above computations do {\it not\1} determine the $\mu'''_{{\rm max},k}$, $\mu'''_{{\rm def},k}$ uniquely. This is closely related to (3.5.3) for $m=2$, which implies that
$$p^{(2)}_k+p^{(2)}_{k-2}-p^{(0)}_k=2\1p^{(2)}_{k-1}\q\q(k\in\Z).
\leqno(3.7.2)$$
Calculating the $\rho_k$ by a computer and applying (3.6.1) in Corollary~(3.6), we can determine the $\mu'''_{{\rm max},k}$ (and hence the $\mu'''_{{\rm def},k}$). If condition~(AT) in the introduction is satisfied, it is enough to calculate the $\rho_k$ for $k\les k_0$ with $\rho_{k_0}-\rho_{k_0-1}=\tau_{Z'}$, where $\tau_{Z'}$ is the Tjurina number of a general hyperplane section $Z'$ of $Z$. For the moment we do not know any examples with $\mu'''_{{\rm def},k}\ne 0$ assuming condition~(AT), see also Remark~(iii) below. This would occur if the 1-dimensional singular locus of $Z\subset\PP^3$ is {\it disconnected}.
\msn
{\bf Caution 2.} There is {\it no method\1} to determine the $\nu'''_{{\rm max},k}$, $\nu'''_{{\rm def},k}$ if the middle morphism $\phi$ in the exact sequence (3.4.5) can be non-vanishing (for instance, in the case of Examples~(5.6--7) below). This problem may be avoided (as in the case of Example~(5.8) below) if the following holds:
$$\mu'_k\,\mu'_{nd-k}=0\q\q(\forall\,k\in\Z),
\leqno(3.7.3)$$
since the middle morphism $\phi$ in (3.4.5) vanishes in this case.
\ms
(iii) In the case condition~(AT) in the introduction does not hold, we may have $\mu'''_{\rm def}\ne 0$, for instance, $\mu'''_{\rm def}=2\1v^6$ in the notation of (4.2) below if $f=x^4z+y^4z+x^2y^2w$, and $\mu'''_{\rm def}=2\1v^7+v^8$ if $f=x^6+x^4yz+y^3w^3$, see also Example~(5.9) below.
\ms
(iv) If $Z\subset\PP^3$ is a free, locally weighted homogeneous divisor (see (4.1) below), we have
$$\mu''_k=\nu''_k=0,
\leqno(3.7.4)$$
and this makes the calculation easier. If $Z$ is {\it strongly free,} the situation is quite simplified, see Remark~(iii) after (4.1) below.
\ms
In the case $n=3$, we can determine the $\rho_k$ quite easily as follows.
\msn
{\bf 3.8.~Proposition.} {\it Let $\ddt=\deg f_{\rm red}$, the {\it reduced degree} of $f$, that is, $\ddt=\sum_j\deg f_j$ if $f=\prod_jf_j^{m_j}$ with $f_j$ irreducible and $m_j\ges1$. Set
$$d_0:=d-\ddt,\q d_1:=2d+\ddt,\q d_2:=3d-1.$$
If $n=3$, then we have
$$\rho_k=\msum_{i\,=\,d_1}^{d_2}\,p^{(2)}_{_k-i}\q\q(k\in\Z),
\leqno(3.8.1)$$
or equivalently}
$$\rho_k=\begin{cases}p^{(3)}_{k-d_1}&\h{\it if}\,\,\,\,k\les d_2,\\ 
p^{(3)}_{d_2-d_1}+d_0(k-d_2)&\h{\it if}\,\,\,\,k>d_2.\end{cases}
\leqno(3.8.2)$$
\msn
{\it Proof.} It is easy to show the equivalence between (3.8.1) and (3.8.2) by using (3.5.2). By Corollary~(3.6), the assertion (3.8.1) is equivalent to the following:
$$\mu'''_{{\max},\,k}=\msum_{i\,=\,3}^{d_0+2}\,p^{(2)}_{k-i}\q\q(k\in\Z).
\leqno(3.8.3)$$
To show the latter equality, it is enough to prove the isomorphism
$$\bl(M'''(3)\br)^{\sim}_{\PP^2}=\OO_{\PP^2}/\OO_{\PP^2}(-Z+Z_{\rm red}).
\leqno(3.8.4)$$
Indeed, this isomorphism implies the following short exact sequences for $k\in\Z$:
$$0\to H^0\bl(\PP^2,\OO_{\PP^2}(k-d_0)\br)\to H^0\bl(\PP^2,\OO_{\PP^2}(k)\br)\to H^0\bl(\PP^2,(M''')^{\sim}_{\PP^2}(k+3)\br)\to 0,
\leqno(3.8.5)$$
where the surjectivity of the last morphism of (3.8.5) follows from the vanishing:
$$H^1\bl(\PP^2,\OO_{\PP^2}(k-d_0)\br)=0\q\q(\forall\,k\in\Z).
\leqno(3.8.6)$$
Hence (3.8.3) follows from (3.5.2) (with $m=3$) and (3.1.9). Here the sheaf $\bl(M'''(3)\br)^{\sim}_{\PP^2}$ on $\PP^2$ is associated with the {\it shifted} graded module $M'''(3)$, or equivalently, with $M'''_{\rm max}(3)$, see \cite[II, Exercise 5.9(c)]{Ha}.
\sk
By the definition of $M'''$ in (3.1), the isomorphism (3.8.4) is reduced to the following vanishing for any finite subset $\Sigma\subset Z$:
$$\Hc^0_{\Sigma}\bl(\OO_{\PP^2}/\OO_{\PP^2}(-Z+Z_{\rm red})\br)=0.
\leqno(3.8.7)$$
(Indeed, (3.8.4) holds outside a finite subset of $Z$ by definition.) This vanishing can be reduced further to the following trivial assertion
$$\Hc^0_{\Sigma}\bl(\OO_{Z_j}(i)\br)=0\q\q(i\in\Z),
\leqno(3.8.8)$$
for any reduced irreducible components $Z_j$ of the non-reduced divisor $Z\subset\PP^2$ by using a finite filtration on $\OO_{\PP^2}/\OO_{\PP^2}(-Z+Z_{\rm red})$. So the assertion (3.8.7) follows. This finishes the proof of Corollary~(3.8). 
\msn
{\bf 3.9.~Remarks.} (i) In the non-reduced case with $n=3$, we can determine the $\mu_k$, $\nu_k$, $\rho_k$ using Proposition~(3.8) and a computer program like Macaulay2 or Singular. In the case of Example~(5.5) below with Macaulay2, this can be done, for instance, as follows:
\ms
\vbox{\small\sf\verb#R=QQ[x,y,z]; f=x^4*y^2*z*(x+y+z)*(x+y);#
\sk
\verb#d=first degree f; I=radical ideal(f); rd=first degree I_0; mxd=3*d;#
\sk
\verb#F=frac(QQ[v]); seq=(1-v^mxd)/(1-v); gam=((v^d-v)/(v-1))^3;#
\sk
\verb#A=QQ[v]/(v^mxd); seqA=sub(seq,A); d1=2*d+rd; d2=3*d;#
\sk
\verb#muA=v^3*sub(hilbertSeries(R/ideal(jacobian ideal(f)),Order=>mxd),vars A)#
\sk
\verb#rhoA=(v^d1-v^d2)*seqA^3#
\sk
\verb#nuA=muA+rhoA-sub(gam,A)#}
\msn
If one copies this from a pdf file, Note in Remark~(3.7)(ii) also applies. (One can replace mxd=3*d by mxd=3*d+3, for instance, if one wants to see more coefficients.) With Singular, one may do, for instance, as follows:
\ms
\vbox{\small\sf\verb#LIB "primdec.lib";#
\sk
\verb#int i, d, rd, mxd, d1, d2; intvec k, mu, nu, rho, gam;#
\sk
\verb#ring R = 0, (x,y,z), dp; poly f=x^4*y^2*z*(x+y+z)*(x+y);#
\sk
\verb#d=deg(f); poly g=x^d+y^d+z^d; ideal J=jacob(f); ideal K=jacob(g);#
\sk
\verb#ideal I=radical(f); rd=deg(I[1]); d1=2*d+rd; d2=3*d; mxd=3*d;#
\sk
\verb#for (i=1; i<mxd-2; i++) {#
\sk\q
\verb#if(i<d1-2) { rho[i]=0; }#
\sk\q
\verb#if (i>=d1-2 && i<d2-2) { rho[i]=(i-d1+3)*(i-d1+4) div 2; }#
\sk\q
\verb#if (i>=d2-2) { rho[i]=(d2-d1)*(d2-d1+1) div 2+(i-d2+3)*(d-rd); }#
\sk\q
\verb#mu[i]=size(kbase(std(J),i-1)); gam[i]=size(kbase(std(K),i-1));#
\sk\q
\verb#k[i]=i+2; nu[i]=mu[i]+rho[i]-gam[i]; }#
\sk
\verb#sprintf(" k : %s",k); sprintf("gam: %s",gam); sprintf("mu : %s",mu);#
\sk
\verb#sprintf("nu : %s",nu); sprintf("rho: %s",rho);#}
\msn
Here Note in Remark~(3.7)(ii) applies.
\ms
(ii) It is very hard to generalize Proposition~(3.8) to the case $n\ges 4$. One difficulty is closely related to the following:
\msn
{\bf Problem.} For a locally free sheaf $L$ on $\PP^1$ having a finite increasing filtration $G$ such that $\Gr^G_iL=\OO_{\PP^1}(a_i)$ with $a_i\in\Z$, this filtration does not necessarily split unless the following splitting condition is satisfied:
$$a_j\les a_i+1\q\h{if}\q i<j.
\leqno{\rm (SC)}$$
\ms
This is closely related to the nonvanishing of $H^1\bl(\PP^1,\OO_{\PP^1}(k)\br)$ for $k\les-2$ (using a spectral sequence), since
$${\rm Ext}^1_{\OO_{\PP^1}}\bl(\OO_{\PP^1}(a_j),\OO_{\PP^1}(a_i)\br)=H^1(\PP^1,\OO_{\PP^1}(a_i-a_j)\br).$$
We have, for instance, the following:
\msn
{\bf Example.} Set $L:=\OO_{\PP^1}(1)\oplus\OO_{\PP^1}(1)$. It has a subsheaf $L_1:=\OO_{\PP^1}\subset L$ generated by the section $(1,x)$. Here $x$ is the affine coordinate of $\C\subset\PP^1$ having a pole of order 1 at $\infty\in\PP^1$, and $\OO_{\PP^1}(1)$ is identified with $\OO_{\PP^1}(\infty)$. We see that $L_2:=L/L_1=\OO_{\PP^1}(2)$, and the following short exact sequence does not split:
$$0\to\OO_{\PP^1}\to\OO_{\PP^1}(1)\oplus\OO_{\PP^1}(1)\to\OO_{\PP^1}(2)\to 0.$$
\sk
Even in the hyperplane arrangement case, it seems quite difficult to construct a filtration on $(M''')^{\sim}_{\PP^{n-1}}$ such that its direct image by $\overline{\pi}_2$ in (3.1.5) gives a filtration satisfying (SC). Even in the case where $n=3$ and $Z$ is a line arrangement having only ordinary double points as singularities, one cannot construct such a filtration {\it on the level of} $(M''')^{\sim}_{\PP^2}$. It is, however, rather interesting that in some cases (for instance, Examples~(5.6) and (5.8) below) we can construct a filtration on $(M''')^{\sim}_{\PP^3}$ such that the graded pieces of a refinement of its direct image by $\overline{\pi}_2$ {\it correctly} give the direct factors although (SC) is not satisfied.
\bs\bs
\vbox{\centerline{\bf 4. Refinements of Theorems~2 and 3}
\bsn
In this section we prove Theorems~(4.3--4) and (4.10) improving Theorems~2--3 in certain cases, and show a symmetry of a modified pole-order spectrum for strongly free divisors in Theorem~(4.5).}
\msn
{\bf 4.1.~Strongly free divisors.} We say that a reduced hypersurface $Z\subset\PP^{n-1}$ is a {\it strongly free\1} divisor, if its {\it affine cone} in $\C^n$ is a free divisor in the sense of \cite{SaK}, that is, the sheaf of vector fields on $\C^n$ which are tangent to the affine cone at smooth points is a free sheaf on $\C^n$ (or equivalently, free at $0\in\C^n$, using the natural $\C^*$-action).
\msn
{\bf Remarks.} (i) A reduced effective divisor $Z\subset\PP^{n-1}$ is strongly free if and only if the kernel of $\df\!\sw\!:\Om^{n-1}\to\Om^n$ is a free graded $R$-module of rank $n-1$ (corresponding to algebraic vector fields on $\C^n$ annihilating $f$). This can be shown by using the Euler vector field $\xi$ as in the proof of Proposition~(4.7) below, see also \cite[Section 8.1]{Di2}, \cite{DiSt2}, etc.
\ms
(ii) Any free divisor on $\PP^{n-1}$ is not necessarily strongly free. For instance, nonsingular hypersurfaces are free, but they are {\it not\1} strongly free if $n\ges 3$, $d\ges 2$. This can be shown by using Remark~(i) above together with the {\it non-vanishing} of the higher extension group ${\rm Ext}_{\OO_X}^n(\OO_X/(\dd f),\OO_X)$, where $(\dd f)\subset\OO_X$ is the Jacobian ideal generated by the partial derivatives of a defining polynomial $f$ of $Z$, which form a {\it regular sequence} of $\OO_X$. It seems well-known that any plane curve is a free divisor (since the sheaf of logarithmic vector fields is reflective) although it is strongly free if and only if $M'=0$ (that is, $M=M''$ in the notation of \cite{DiSa2}) by an argument similar to Remark~(iii) below, see also \cite{DiSt2}, etc. As an example with $n=4$, we have, for instance, a {\it generic\1} hyperplane arrangement in $\PP^3$ with degree $5$, which is free, but not strongly free. (Here ``generic" means a divisor with {\it normal crossings.}) These give examples of {\it non-splitting\1} vector bundles on $\PP^{n-1}$, see for instance \cite{DiSt2}, etc.
\ms
(iii) If $Z$ is a strongly free reduced divisor on $\PP^3$, then the graded $R$-module $M$ in (3.1) is a 2-dimensional {\it Cohen-Macaulay} $R$-module by Remark~(i) after (4.1) above (using Auslander--Buchsbaum formula, see \cite{BrHe}, \cite[Theorem 19.9]{Ei1}). Here $M$ is isomorphic to the graded $R$-algebra $R/(\dd f)$ as graded $R$-module (up to a shift of grading), and the direct image of the associated sheaf $M^{\sim}$ by the projection $\pi_2$ in (3.1) is free at the origin (using for instance \cite[II, Theorem 8.21A(c)]{Ha}, \cite[Theorem 6.8]{Ei1}). This implies that
$$M'=M''=M'''_{\rm def}=N'=N''=N'''_{\rm def}=0.
\leqno(4.1.1)$$
Indeed, the assertions for $M',M'',M'''_{\rm def}$ follow from the above freeness of $((\pi_2)_*M^{\sim})_0$, and these imply the remaining by Theorem~(3.4). We then get the isomorphisms
$$M=M'''_{\rm max},\q N=N'''_{\rm max},\q Q=Q'''_{\rm max},
\leqno(4.1.2)$$
where the last assertion for $Q$ is by Theorem~(3.4).
\ms
(iv) We have a {\it symmetry\1} of $\R_f$ with center 1 if $Z$ is a strongly free, locally positively weighted homogeneous divisor, see \cite{Na}. In this paper we say that a reduced divisor $Z$ on a complex manifold $Y$ is {\it locally positively weighted homogeneous,} if $(Z,z)\subset(Y,z)$ is defined by a weighted homogeneous polynomial with {\it strictly positive} weights using some local coordinate system of $(Y,z)$ at {\it any} point $z\in Z$. This is also called {\it strongly} (or {\it locally}) {\it quasi-homogeneous,} see \cite{CNM}, \cite{DiSt2}.
\msn
{\bf 4.2.~Case $n=4$.} Assume $n=4$. In the notation of (3.5), consider the formal power series
$$\mu:=\msum_k\,\mu_k\1v^k,\q\nu:=\msum_k\,\nu_k\1v^k,\q\rho:=\msum_k\,\rho_k\1v^k\q\h{in}\,\,\,\Z[[v]],$$
and similarly for $\mu'$, $\mu''$, $\mu'''_{\rm def}$, $p^{(m)}$, etc. Define
$${\rm Diff}^2\mu:=(1-v)^2\mu\q\h{in}\,\,\,\Z[[v]],$$
so that
$$\mu=(1-v)^{-2}\1{\rm Diff}^2\mu\q\h{in}\,\,\,\,\Z(\!(v)\!)\,\bl(:=\Z[[v]][v^{-1}]\br),$$
(similarly for $\nu$, $\rho$). Note that Diff in this section corresponds to the one in (3.5), and
$$\aligned&(1\,{-}\,v)^{-1}=1\,{+}\,v\,{+}\,v^2{+}\cdots\q\h{in}\,\,\,\Z(\!(v)\!),\\&p^{(m)}=(1\,{-}\,v)^{-m}=\bl(p^{(1)}\br)^m\,\,\,\,\,(m\ges 0).\endaligned$$
\sk
If $Z$ is a {\it strongly free\1} divisor, then (4.1.2) together with (3.1.7) implies that
$${\rm Diff}^2\mu,\,\,{\rm Diff}^2\nu,\,\,{\rm Diff}^2\rho\,\,\,\h{are polynomials in $v$ with {\it positive} coefficients.}
\leqno(4.2.3)$$
The {\it positivity\1} of coefficients does not necessarily holds unless $Z$ is strongly free; for instance, if $f=xyzw(x+y+z)(y-z+w)$ in \cite[Example 5.7]{DiSt2}.
\sk
In the {\it strongly free\1} divisor case, we get the following by the duality isomorphisms in Theorem~(3.4) (or more concretely, by (3.6.1) in Corollary~(3.6)) together with (4.1.2):
$${\rm Diff}^2\rho(v)=v^{4d+2}({\rm Diff}^2\mu)(v^{-1})\q\h{in}\,\,\,\Z(v),
\leqno(4.2.4)$$
where the multiplication by $v^2$ comes from the dimension of (the support of) $M$.
\sk
This determines $\nu$ using the following well-known relation in $\Z[[v]]$\,:
$$\mu-\nu+\rho=\gamma:=\bl(\tfrac{v^d-v}{v-1}\br)^4\,\in\,\Z[v].
\leqno(4.2.5)$$
This relation follows from the assertion that the Euler characteristic of a bounded complex of finite dimensional vector spaces is {\it independent of its differential.} Note that $\gamma$ calculates the top cohomology of the Koszul complex for a nonsingular hypersurface $Z_{\rm ns}\subset\PP^3$, and gives the {\it spectrum\1} of the isolated hypersurface singularity of the cone of $Z_{\rm ns}$ in $\C^4$.
\sk
The {\it Euler characteristic} of the $E_1$-complex of the pole order spectral sequence can be defined by the formal power series
$$\chi_f:=\mu-v^{-d}\1\nu+v^{-2d}\1\rho\,\,\in\,\,\Z[[v]].
\leqno(4.2.6)$$
Here we have the divisions by $v^d$ and $v^{2d}$, since the $E_1$-differential $\ddd^{(1)}$ decreases the degree by $-d$. Setting $\mu^{(2)}=\msum_k\,\mu^{(2)}_k\,v^k$, etc., we have the equality
$$\chi_f=\mu^{(2)}-v^{-d}\nu^{(2)}+v^{-2d}\rho^{(2)}\q\h{in}\,\,\,\Z[[v]].
\leqno(4.2.7)$$
So $\chi_f$ is also called the {\it modified pole-order spectrum for the $E_2$-term}. Note that
$$\chi_f\in\Z[v],\q\h{that is,}\q\mu_{k-d}+\rho_{k+d}=\nu_k\,\,\,\,\h{for}\,\,\,k\gg 0.
\leqno(4.2.8)$$
Indeed, $\mu_k+\rho_k=\nu_k$ and $\mu_{k+1}=\mu_k+r'''$ (similarly for $\rho_k$) for $k\gg 0$ by (4.2.5), (3.5.5).
\sk
For positive integers $m$, set
$$\aligned{\rm Eu}_{\1i}^{\les m}&:=\bl(\msum_{j\in\N,\,i+jd\les m}\,\chi_{f,i+jd}\br)-\delta_{i,d},\\{\rm Eu}_{\1i}&:={\rm Eu}_{\1i}^{\les m}\q\h{for}\q m\gg 0,\endaligned
\leqno(4.2.9)$$
where $\chi_f=\msum_k\,\chi_{f,k}\,v^k$.
In the notation of Theorem~3, we have
$${\rm Eu}_{\1i}^{\les 2d-2}=\cho^{(2)}_{f,\,\ee(-i/d)}\bl(\twd\br).
\leqno(4.2.10)$$
\sk
Theorem~(4.3) below with $m=2d-2$ improves Theorem~3 for hyperplane arrangements or strongly free, locally positively weighted homogeneous divisors on $\PP^3$. In the latter case, it is known that condition~(4.3.1) is satisfied (see Remark~(ii) after the proof of Theorem~(4.3) below), and moreover conditions~(4.3.2--3) can be proved, see Corollary~(4.6), Theorem~(4.10) below. (Note that hyperplane arrangements are locally positively weighted homogeneous.)
\msn
{\bf 4.3.~Theorem.} {\it Let $Z$ be a reduced divisor on $\PP^3$ and $m\in\Z$ satisfying conditions~{\rm (GH)} and $m\in[2d-2,2d+3]$. Assume the following conditions hold in the notation of $(4.2)\,{:}$
$$\max\R_f\les\tfrac{m}{d},
\leqno(4.3.1)$$
$$\deg\chi_f\les m,
\leqno(4.3.2)$$
$${\rm Eu}_{\1i}=-\chi(U)\q(i\in[1,d]),
\leqno(4.3.3)$$
Then the pole order spectral sequence degenerates completely at $E_3$, and almost at $E_2$, more precisely, we have}
$$\mu_k^{(2)}=\nu_{k+d}^{(2)}=\rho_{k+2d}^{(2)}=0\q\h{if}\q k>m.
\leqno(4.3.4)$$
\msn
{\it Proof.} We have
$$\mu_k^{(\infty)}=\nu_{k+d}^{(\infty)}=\rho_{k+2d}^{(2)}=0\q\h{if}\q k>m.
\leqno(4.3.5)$$
Indeed, the assertions for $\mu_k^{(\infty)}$, $\rho_{k+2d}^{(2)}$ follow from the condition~(4.3.1) (using Theorem~1) and Proposition~1. The assertion for $\nu_{k+d}^{(\infty)}$ is then proved by using the inequalities 
$$-\chi(U)\les{}^{(\infty)}{\rm Eu}_{\1i}^{\les m}\les{\rm Eu}_{\1i}^{\les m}={\rm Eu}_{\1i}=-\chi(U),
\leqno(4.3.6)$$
which is a generalization of (2.6.1). Here $^{(\infty)}{\rm Eu}_{\1i}^{\les m}$ is defined by replacing $\chi_f$ in (4.2.9) with $\chi_f^{(\infty)}$, and the latter is defined by replacing the $\mu^{(2)}_k$, etc.\ in (4.2.7) with $\mu^{(\infty)}_k$, etc. The first two inequalities can be shown as in (2.6). The last two equalities follow from conditions~(4.3.2--3). Note that the first inequality of (4.3.6) becomes a strict inequality if $\nu_{k+d}^{(\infty)}\ne 0$ for some $k>m$.
\sk
We now show (4.3.4). The assertion for $\rho_{k+2d}^{(2)}$ is proved in (4.3.5). We have the equalities
$$\mu^{(2)}_k=\nu^{(2)}_{k+d}\q(k>m),$$
using the equality (4.2.7) and the vanishing of $\chi_{f,k}$ coming from (4.3.2) together with the assertion for $\rho_{k+2d}^{(2)}$ in (4.3.5). If $\mu^{(2)}_k\ne 0$ for some $k>m$, we can take the minimal $k>m$ satisfying
$$\mu^{(2)}_k=\nu^{(2)}_{k+d}\ne 0.$$
Since $\nu_{k+d}^{(\infty)}=0$ by (4.3.5), we then get the non-triviality of a differential
$$\ddd^{(r)}:N^{(r)}_{k+d}\to M^{(r)}_{k-(r-1)d}\q(\exists\,r\ges 2).$$
Here $k':=k-(r-1)d<m$ by the above assumption on $k>m$ (since $\mu_{k'}
^{(r)}\ne 0$). This deduces a contradiction using the inequalities in (4.3.6) together with (4.3.3). Indeed, the second inequality of (4.3.6) becomes a strict inequality in this case. So (4.3.4) follows.
\sk
The $E_3$-degeneration now follows from (4.3.4), since $\ddd^{(r)}$ decreases the degree by $rd$ (and $\nu^{(r)}_k$ vanishes for $k\les d$ by definition). This finishes the proof of Theorem~(4.3).
\msn
{\bf 4.4.~Remarks.} (i) The $E_3$-degeneration in Theorem~(4.3) holds for {\it any} degree, and not only for some {\it lower\1} degrees as in Theorem~2 or 3. If conditions~(GH) and (4.3.1--2) are satisfied, then condition~(4.3.3) is equivalent to condition~(4.3.4), since the latter implies that
$${}^{(\infty)}{\rm Eu}_{\1i}^{\les m}={\rm Eu}_{\1i}^{\les m}\q(i\in[1,d]).$$
In the case conditions~(4.3.3--4) are not satisfied, the spectral sequence does not degenerate at $E_2$ even after restricting to any sufficiently high degrees.
\ms
(ii) Condition~(4.3.1) for $m=2d-2$ is satisfied in the case $Z$ is a hyperplane arrangement (see \cite[Theorem~1]{bha}) or a strongly free, locally positively weighted homogeneous divisor (using Remark~(iv) after (4.1)).
\ms
(iii) In the case of strongly free, locally positively weighted homogeneous divisors on $\PP^{n-1}$, it will be shown that $\chi_{f,d}=-\chi(U)+1$ in Corollary~(4.8) below. So (4.3.3) can be reduced to the independence of $i$ of the ${\rm Eu}_{\1i}$ if condition~(4.3.2) with $m<2d$ is satisfied.
\ms
(iv) By Corollary~(3.6) and Remark~(3.7)(ii), Theorem~(4.3) also applies to hyperplane arrangements which are not necessarily strongly free, if we can compute the $\rho_k$ for $k\les k_0$ with $\rho_{k_0}-\rho_{k_0-1}=\tau_{Z'}$, where $\tau_{Z'}$ is the Tjurina number of a general hyperplane section $Z'$ of $Z$, see for instance Example~(5.8) below.
\ms
(v) It is not easy to generalize Theorem~(4.3) to the generic $E_3$-degeneration case, since it is unclear whether the $\mu^{(2)}_k$ are weakly decreasing for $k\ges 3d-2$.
\ms
As a corollary of Theorem~(3.4), we get the following (which may be related to Remark~(iv) after (4.1)).
\msn
{\bf 4.5.~Theorem.} {\it If $Z$ is a strongly free reduced divisor on $\PP^3$, we have the symmetry of the modified pole-order spectrum for the $E_2$-term}
$$\chi_f(v)=v^{2d}\1\chi_f(v^{-1})\q\h{in}\,\,\,\Q(v).
\leqno(4.5.1)$$
\msn
{\it Proof.} Set
$$\chit_f:=(1-v)^2\chi_f,\q\mut:=(1-v)^2\mu,\q\nut:=(1-v)^2\nu,\q\rhot:=(1-v)^2\rho.$$
These belong to $\Z[v]$ by (4.1.2) and (3.1.7). Since
$$v^2(1-v^{-1})^2=(1-v)^2\q\h{in}\,\,\,\Q(v),$$
the assertion (4.5.1) is equivalent to
$$v^{2d+2}\1\chit_f(v^{-1})=\chit_f(v)\q\h{in}\,\,\,\Q(v),
\leqno(4.5.2)$$
\sk
By Theorem~(3.4) (or rather Corollary~(3.6)) together with (4.1.2), we get
$$v^{4d+2}\1\mut(v^{-1})=\rhot(v),\q v^{4d+2}\1\nut(v^{-1})=\nut(v).
\leqno(4.5.3)$$
These imply that
$$v^{2d+2}\1\mut(v^{-1})=v^{-2d}\1\rhot(v),\q v^{2d+2}\bl(v^d\1\nut(v^{-1})\br)=v^{-d}\1\nut(v).
\leqno(4.5.4)$$
So (4.5.2) follows. This finishes the proof of Theorem~(4.5).
\ms
Theorem~(4.5) implies the following.
\msn
{\bf 4.6.~Corollary.} {\it If $Z$ is a strongly free reduced divisor on $\PP^3$, then {\rm (4.3.2)} holds.}
\msn
{\it Proof.} It is enough to consider the case $m=2d-2$ (since we assume $m\ges 2d-2$). By the symmetry in Theorem~(4.5), the assertion (4.3.2) is reduced to
$$\chi_{f,k}=0\q\h{if}\,\,\,\,k\les 1.
\leqno(4.6.1)$$
By the definition in the introduction, however, we have
$$M_k=0\q(k\les 3),\q N_{k+d}=0\q(k\les 2),\q Q_{k+2d}=0\q(k\les 1),
\leqno(4.6.2)$$
since $\Om^p_k=0$ for $k<p$ in general. So the assertion follows from the definition of $\chi_f$ in (4.2.6). This finishes the proof of Corollary~(4.6).
\ms
The following is closely related to condition~(4.3.3) for $i=d$, see Corollary~(4.8) below.
\msn
{\bf 4.7.~Proposition.} {\it In the notation of the introduction, let $(A_f^{\ssb},\ddd)\subset(\Om^{\ssb},\ddd)$ be the graded subcomplex with $A_f^p:={\rm Ker}(\df\!\sw\!:\Om^p\to\Om^{p+1})$. Assume $Z$ is a strongly free, locally positively weighted homogeneous divisor on $\PP^{n-1}$ or a hyperplane arrangement in $\PP^3$. Then we have
$$\dim H^j(A_{f,d}^{\ssb},\ddd)=\dim H^{j-1}(U,\C)\q(j\in\Z),
\leqno(4.7.1)$$
with $(A_{f,d}^{\ssb},\ddd)$ the degree $d$ part of the graded complex $(A_f^{\ssb},\ddd)=\mopl_{k\in\Z}\,(A_{f,k}^{\ssb},\ddd)$.}
\msn
{\it Proof.} Set $D:=f^{-1}(0)\subset X:=\C^n$. Let $(\Om^{\ssb}(\log D),\ddd)$ be the logarithmic complex defined by $f\om\in\Om^p$, $f\ddd\om\in\Om^{p+1}$, where the second condition is equivalent to $\df\!\sw\om\in\Om^{p+1}$ assuming the first, see \cite{SaK}. We have the direct sum decomposition
$$\Om^p(\log D)=\onf\1A_f^p\oplus\onf\1\ix\1A_f^{p+1}\q(p\in\Z),
\leqno(4.7.2)$$
where $\ix$ denotes the inner derivation by the Euler field
$$\xi:=\msum_{i=1}^n\,\ond\1x_i\tfrac{\dd}{\dd x_i}.$$This follows from the well-known relations
$$\ix(\df\!\sw\om)+\df\!\sw\ix\1\om=f\om\q\bl(\om\in\Om^p\bl[\onf\br]\br),
\leqno(4.7.3)$$
$$\ix(\ddd\om)+\ddd(\ix\1\om)=\kod\,\1\om\q\bl(w\in\Om^p\bl[\onf\br]_k\br).
\leqno(4.7.4)$$
Indeed, if $\df\!\sw\eta\in f\1\Om^{p+1}$ for $\eta:=h\1\om\in\Om^p$, then (4.7.3) implies that
$$\df\!\sw\ix\bl(\onf\df\!\sw\eta\br)=\df\!\sw\eta,\q\h{hence}\q \ix\bl(\onf\df\!\sw\eta\br)-\eta\in A_f^p.$$
Note that $\ix\1A_f^{p+1}=A_f^{p+1}$ in (4.7.2), since
$$\bl(\onf\1\df\!\sw\br)\ssc\ix={\rm id}\q\h{on}\q A_f^{p+1}.
\leqno(4.7.5)$$
So (4.7.4) gives an isomorphism of complexes
$$(\Om^{\ssb}(\log D)_k,\ddd)=C\bl(\kod:(A_{f,k+d}^{\ssb},\ddd)\to(A_{f,k+d}^{\ssb},\ddd)\br)\q(k\in\Z),
\leqno(4.7.6)$$
with left-hand side the degree $k$ part of the graded complex $(\Om^{\ssb}(\log D),\ddd)$. (Note that $\onf$ in (4.7.2) does not affect the differential $\ddd$, since $A_f^p={\rm Ker}\,\df\!\sw$.) We thus get the quasi-isomorphism
$$(\Om^{\ssb}(\log D),\ddd)\cong(A_{f,d}^{\ssb},\ddd)[1]\oplus(A_{f,d}^{\ssb},\ddd),
\leqno(4.7.7)$$
inducing the isomorphisms
$$H^j(\Om^{\ssb}(\log D),\ddd)=H^{j+1}(A_{f,d}^{\ssb},\ddd)\oplus H^j(A_{f,d}^{\ssb},\ddd)\q(j\in\Z).
\leqno(4.7.8)$$
By \cite[Theorem 1]{CNM} and \cite[Corollary 6.3]{WiYu}, we have the isomorphisms
$$H^j(\Om^{\ssb}(\log D),\ddd)=H^j(X\setminus D,\C)\q(j\in\Z),
\leqno(4.7.9)$$
(replacing the direct sum with the convergent infinite sum if necessary). We have moreover the isomorphisms
$$H^j(X\setminus D,\C)=H^j(U,\C)\oplus H^{j-1}(U,\C)\q(j\in\Z),
\leqno(4.7.10)$$
since $X\setminus D$ is a $\C^*$-bundle over $U$ and its Chern class with $\Q$-coefficients vanishes. So (4.7.1) follows by decreasing induction on $j$. This finishes the proof of Proposition~(4.7).
\msn
{\bf Remark.} The above argument is a generalization of an argument in the proof of \cite[Theorem 3.9]{DiSt2} (which gives an assertion only for the highest cohomology).
\msn
{\bf 4.8.~Corollary.} {\it In the notation of $(4.2.9)$ and under the assumption of Proposition~$(4.7)$ with $n=4$, we have}
$$\chi_{f,d}=-\chi(U)+1.
\leqno(4.8.1)$$
\msn
{\it Proof.} Since $\df\!\sw\eta$ ($\eta\in\Om^p$) has degree at least $d+p$, we have
$$H^j(\Om^{\ssb},\df\!\sw)_d=A_{f.d}^j\q(j>0),$$
So the assertion follows from Proposition~(4.7).
\ms
The following is a key to the proof of Theorem~(4.10) below.
\msn
{\bf 4.9.~Proposition.} {\it For each $i\in[1,d]$, the following two conditions are equivalent\,$:$
$${\rm Eu}_{\1i}=-\chi(U),
\leqno(4.9.1)$$
$$\mu_{\1i+jd}-\rho_{\1i+jd}-r'''d=\chi(U)-\chi(U_{\rm ns})\q\h{for}\,\,\,j\gg 0,
\leqno(4.9.2)$$
where $r'''$ is as in $(3.5.5)$, and $U_{\rm ns}:=\PP^3\setminus\{f_{\rm ns}=0\}$ with $f_{\rm ns}:=\sum_{i=1}^4\,x_i^d$.}
\msn
{\it Proof.} For $g=\msum_k\,g_k\1v^k\in\Q[[v]]$ and $m\gg 0$, set
$${\rm Eu}_{\1i}^{\les m}(g)=\msum_{j\in\Z,\,i+jd\les m}\,g_{i+jd},$$
For {\it any\1} sufficiently large integer $m$, we have by (4.2.8--9)
$${\rm Eu}_{\1i}={\rm Eu}_{\1i}^{\les m-2d}(\mu)-{\rm Eu}_{\1i}^{\les m-d}(\nu)+{\rm Eu}_{\1i}^{\les m}(\rho)-\delta_{i,d}\q(i\in[1,d]).
\leqno(4.9.3)$$
Note that (2.4.1) holds also for the isolated singularity $f_{\rm ns}$, $U_{\rm ns}$, and moreover $\gamma$ in (4.2.5) coincides with the spectrum $\chi_{f_{\rm ns}}$ of the isolated singularity $f_{\rm ns}$, where (4.3.3) holds for $f_{\rm ns}$, $U_{\rm ns}$. For any $k\gg 0$, we then get by (4.2.5)
$$-\chi(U_{\rm ns})={\rm Eu}_{\1i}^{\les m}(\mu)-{\rm Eu}_{\1i}^{\les m}(\nu)+{\rm Eu}_{\1i}^{\les m}(\rho)-\delta_{i,d}\q(i\in[1,d]).
\leqno(4.9.4)$$
So the assertion follows by comparing (4.9.3--4). Indeed, the difference is given by
$$\mu_k+\mu_{k-d}-\nu_k=\mu_k-\rho_k-r'''d\q\h{for some}\,\,\,k=i+jd\gg 0.$$
This finishes the proof of Proposition~(4.9.1).
\ms
We now get the following (which can be applied to \cite[Examples 5.3--5]{DiSt2}, for instance).
\msn
{\bf 4.10.~Theorem.} {\it Let $Z$ be a strongly free, locally positively weighted homogeneous divisor on $\PP^3$. Then the pole order spectral sequence degenerates completely at $E_3$, and almost at $E_2$, more precisely, $(4.3.4)$ holds for $m=2d-2$.}
\msn
{\it Proof.} This follows from Theorem~(4.3) for $m=2d-2$, Corollaries~(4.6), (4.8), and also Proposition~(4.9). Indeed, condition~(4.9.1) holds for $i=d$ by Corollary~(4.8) together with condition~(4.3.2) for $m=2d-2$. So it holds for any $i\in[1,d]$, since $\mu_k-\rho_k$ is {\it independent\1} of $k\gg 0$, see a remark after (4.2.8). This finishes the proof of Theorem~(4.10).
\msn
{\bf 4.11.~Remarks.} (i) In the case of essential hyperplane arrangements in $\PP^3$, we can show that Theorem~(4.10) holds for $m=2d-1$, using Castelnuovo-Mumford regularity. Here it is not easy to control $\mu'''_{\rm def}$. (This subject will be treated in another paper.)
\ms
(ii) It may be possible to deduce the equality (4.8.1) by combining some arguments in the proofs of \cite[Theorem 3.9]{DiSt2} and Theorem~(4.3) in this paper. This is left to the reader.
\ms
(iii) In \cite{DiSt2}, an algorithm calculating the $\mu_k^{(2)}$ for $k<2d$ (with $n=4$) is given. By (4.2.7), we have
$$\mu^{(2)}-\chi_f=v^{-d}\,\nu^{(2)}-v^{-2d}\,\rho^{(2)}.
\leqno(4.11.1)$$
Sometimes symmetry holds for $\mu^{(2)}$, but sometimes not (see {\it loc.\,cit.}) although it always hold for $\chi_f$ by Theorem~(4.5). This depends on $v^{-d}\,\nu^{(2)}$, since it always holds for $v^{-2d}\,\rho^{(2)}$ by \cite[Proposition 2.2]{DiSt2}.
\ms
(iv) For \cite[Examples~5.3--5]{DiSt2}, we can verify that condition~(4.3.3) in Theorem~(4.3) is satisfied by using a computer program like Macaulay2 or Singular. In the case of Example~5.5 in \h{\it loc.\,cit.} with Macaulay2, for instance, we can compute the dimensions of the $E_1$-terms of the pole order spectral sequence very rapidly by typing M2 and pressing RETURN on a terminal of Unix (or Mac, \h{etc.}) with Macaulay2 installed, and then copying and pasting the following (after some corrections if necessary, see Note in Remark~(3.7)(ii)):
\ms
\vbox{\small\sf\verb#R=QQ[x,y,z,w];#
\sk
\verb#f=(x^3-y^3)*(x^3-z^3)*(x^3-w^3)*(y^3-z^3)*(y^3-w^3)*(z^3-w^3);#
\sk
\verb#d=first degree f; mxd=4*d; A=QQ[v]/(v^mxd); F=frac(QQ[v]); vA=sub(v,A);#
\sk
\verb#muA=vA^4*sub(hilbertSeries(R/ideal(jacobian ideal(f)),Order=>mxd),vars A);#
\sk
\verb#seqA=sub((1-v^mxd)/(1-v),A); DDmu=sub(muA*(1-vA)^2,F);#
\sk
\verb#DDrho=sub(DDmu,{v=>1/v})*v^(4*d+2); rhoA=sub(DDrho,A)*seqA^2;#
\sk
\verb#nuA=muA+rhoA-sub(((v^d-v)/(v-1))^4,A); DDnu=sub(nuA*(1-vA)^2,F);#
\sk
\verb#rhoSSA=sub(DDrho/v^(2*d),A)*seqA^2; nuSA=sub(DDnu/v^d,A)*seqA^2;#
\sk
\verb#ChA=muA-nuSA+rhoSSA; Ch=sub(ChA,F)#
\sk
\verb#B=QQ[v]/(v^(d+1)); ChB1=sub(ChA,B); ChA1=sub(ChB1,A);#
\sk
\verb#ChA2=(ChA-ChA1)/vA^d; Eu=ChA1+ChA2-vA^d#
\sk
\verb#EuB=sub(Eu,B); e=(EuB-EuB*v/v)/v^d; Eu-sub(sub(vA*seqA,B)*e,A)#}
\msn
If one copies this from a pdf file, then Note in Remark~(3.7)(ii) applies. Here {\sf Ch} is the {\it modified pole-order spectrum $\chi_f$ for the $E_2$-term} in (4.2.6) which is computed by using (4.2.4--5). If the final output is 0, then this means that condition~(4.3.3) is satisfied. (One can make a similar computation using Singular as in Remark~(3.7)(ii) and (3.9)(i). This is left to the reader.)
\sk
It seems interesting to compare the above {\sf Ch} with the {\it pole order spectrum} $\mu^{(2)}$ for $H^3(F_{\!f})$ at the $E_2$-level, which is calculated in \cite{DiSt2} (where $F_{\!f}$ is the Milnor fiber and $v=t^{1/d}$). By (4.11.1) the difference $\mu^{(2)}\,{-}\,{\sf Ch}$ coincides with $v^{-d}\1\nu^{(2)}-v^{-2d}\1\rho^{(2)}$, that is, the difference between the pole order spectra for $H^2(F_{\!f})$ and $H^1(F_{\!f})$ at the $E_2$-level which are shifted by $d$ and $2d$ respectively. If this is a polynomial of degree at most $3d/2$, then this implies that the pole order spectral sequence degenerates at $E_2$ except for degree $d/2$, using complex conjugation as in (9) in the introduction.
\ms
(v) If conditions (GH) and (AT) hold, for instance, if $Z$ is a hyperplane arrangement in $\PP^3$, we may add the following to the code of Macaulay2 in Remark~(3.7)(ii):
\ms
\vbox{\small\sf\verb#F=frac(QQ[v]); vA=sub(v,A); DmusA=musA*(1-vA)#
\sk
\verb#seqA=sub((1-v^mxd)/(1-v),A); DDmutA=mutA*(1-vA)^2#
\sk
\verb#DDrhoA=sub(sub(sub(DDmutA,F),v=>1/v)*v^(4*d+2),A)#
\sk
\verb#DnusA=sub(sub(sub(DmusA,F),v=>1/v)*v^(4*d+1),A)#
\sk
\verb#nutdA=sub(sub(sub(mupA,F),v=>1/v)*v^(4*d),A)#
\sk
\verb#gamA=sub(((v^d-v)/(v-1))^4,A); rhoA=DDrhoA*seqA^2;#
\sk
\verb#DDnutmA=(muA+rhoA-gamA-DnusA*seqA+nutdA)*(1-vA)^2#
\sk
\verb#nusSA=(DnusA/vA^d)*seqA; nutmSA=(DDnutmA/vA^d)*seqA^2;#
\sk
\verb#nutdSA=nutdA/vA^d; nuSA=nusSA+nutmSA-nutdSA;#
\sk
\verb#rhoSSA=(DDrhoA/vA^(2*d))*seqA^2;#}
\msn
and also the {\it last four lines} of the code in Remark~(iii) above. In the case the last output is zero and $\deg{\sf Ch}<2d$, we may get the degeneration at $E_3$ (almost at $E_2$) by Theorem~(4.3) (choosing $m$ so that (4.3.1--2) holds, where Theorem~2 or 3 may be needed to determine $\R_f$). Here we have to calculate also $\chi(U)$ (except for the hyperplane arrangement case using Corollary~(4.8)), and we {\it assume\1} that $\mu'''_{\rm def}$ {\it vanish}, see also Remark~(3.7)(iii). The last vanishing must be verified by computing the $\rho_k$ for $k\les k_0$ with $\rho_{k_0}-\rho_{k_0-1}=\tau_{Z'}$. As for the condition: $\deg{\sf Ch}<2d$, this may be relaxed if the definition of {\sf Eu} is modified (and $\chi(U)$ can be calculated). Note also that the $\nu'''_{\rm max}$, $\nu'''_{\rm def}$ cannot be calculated correctly by the above code in general, although this does not matter for the computation of {\sf nuSA}.
\bs\bs
\vbox{\centerline{\bf 5. Examples}
\bsn
In this section we calculate some examples.}
\msn
{\bf 5.1.~Example.} Let
$$f=x^5+y^4z+x^3y^2\q\h{with}\q h=x^5+y^4+x^3y^2,$$
where $n=3$, $d=5$. Using computer programs based on \cite{DiSt1} or this paper, we can get
$$\scalebox{0.9}{$\begin{array}{rccccccccccccccc}
k:&3 &4 &5 &6 &7 &8 &9 &10 &11 &12 &13 &14 &\cdots\\
\gamma_k:&1 &3 &6 &10 &12 &12 &10 &6 &3 &1 & &\\
\mu_k:&1 &3 &6 &10 &12 &12 &11 &11 &11 &11 &11 &11 &\cdots\\
\mu^{\scriptscriptstyle(2)}_k:&1 &2 &1 &2 &2 &2 &1 &1 &1 &1 &1 &1 &\cdots\\
\mu^{\scriptscriptstyle(3)}_k:& &1 & &1 &1 &1 & & & & & &\\
\nu_k:& & & & & & &1 &5 &8 &10 &11 &11 &\cdots\\
\nu^{\scriptscriptstyle(2)}_k:& & & & & & & & & & &1 &1 &\cdots\\
\nu^{\scriptscriptstyle(3)}_k:& & & & & & & & & & & &\end{array}$}$$
Here 0 is omitted to simplify the display. This answers a question in \cite[Section 2.4]{ex}. The spectral sequence degenerates at $E_3$ for $k\les 3d-1$ by Theorem~2, and a variant of condition~(18) with $\R_Z$ replaced by $\{d\}$ holds in this case, where $\chi(U)=1$. Since $\R_Z\cap\tfrac{1}{5}\,\Z=\{1\}$, we get by using Corollary~1
$$5\,\R_f^0=\{4,6,7,8\},$$
as in \h{\it loc.~cit.} (where a computer calculation by RISA/ASIR is used).
\msn
{\bf 5.2.~Example.} Let
$$f=(yz-x^2)^3+y^6\q\h{with}\q h=(y-x^2)^3+y^6,$$
where $n=3$, $d=6$. In this case $Z$ has only one singular point which is a $\mu$-constant deformation of a weighted homogeneous polynomial $u^3+y^{12}$ with $\mu_Z=22$, $\tau_Z=21$, see Remark below. Topologically $Z$ is a union of three spheres intersecting at one point so that $\chi(Z)=2\cdot 3-2=4$, and hence $\chi(U)=3-4=-1\,(=5\cdot 4+1-\mu_Z)$ (see \cite[2.9.5]{wh} for the last equality). It is rather rare that we have $\chi(U)<0$ with $n=3$.
\sk
Running computer programs based on \cite{DiSt1} or this paper, we can get
$$\scalebox{0.9}{$\begin{array}{rccccccccccccccccccc}
k:&3 &4 &5 &6 &7 &8 &9 &10 &11 &12 &13 &14 &15 &16 &17 &\cdots\\
\gamma_k:&1 &3 &6 &10 &15 &18 &19 &18 &15 &10 &6 &3 &1 & &\\
\mu_k:&1 &3 &6 &10 &15 &18 &20 &21 &21 &21 &21 &21 &21 &21 &21 &\cdots\\
\mu^{\scriptscriptstyle(2)}_k:&1 &1 &2 &1 &2 &1 &1 &1 &1 &1 &1 &1 &1 &1 &1 &\cdots\\
\mu^{\scriptscriptstyle(3)}_k:&1 & &1 & &1 & & & & & & & & & &\\
\nu_k:& & & & & & &1 &3 &6 &11 &15 &18 &20 &21 &21 &\cdots\\
\nu^{\scriptscriptstyle(2)}_k:& & & & & & &1 &1 &2 &2 &2 &1 &1 &1 &1 &\cdots\\
\nu^{\scriptscriptstyle(3)}_k:& & & & & & &1 &1 &2 &2 &2 &1 &1 & &\end{array}$}$$
The spectral sequence then degenerates at $E_3$ for $k\les 3d-1$ by Theorem~2, and a variant of condition~(18) with $\R_Z$ replaced by $\{d\}$ does not hold in this case.
\sk
On the other hand, a calculation by RISA/ASIR shows
$$12\,\R_Z=\{5,\dots,18\}.$$
This imply that ${\rm CS}(f)=\emptyset$. So we get $\R_f^0=\emptyset$ as a consequence of the above calculation using Corollary~1. This is compatible with a computation using the program ``bfunction" in RISA/ASIR.
\msn
{\bf Remark.} We can calculate the (global) Tjurina number of $h$ using Singular as follows:
\ms
\vbox{\small\sf\verb#ring R = 0, (x,y), dp;#
\sk
\verb#poly h=(y-x^2)^3+y^6;#
\sk
\verb#ideal J=(jacob(h),h);#
\sk
\verb#vdim(groebner(J));#}
\msn
This coincides with the global Tjurina number $\tau_Z$ of $Z$ (that is, $\mu_k,\nu_k$ for $k\gg 0$) if $z$ is chosen so that $\{z=0\}\cap{\rm Sing}\,Z=\emptyset$, see \cite{DiSa2}. Here one can also use {\sf tjurina(h);} after writing {\sf LIB "sing.lib";} at the beginning. If the local Tjurina number at 0 is needed, one must replace {\sf dp} with {\sf ds} (so that a computation is done in the localization of the polynomial ring at the origin); for instance, when {\sf\verb#h=(x^3-y^2)*(x^2-y^3)#}.
We can also get the global Milnor number (in order to see whether the singularities are weighted homogeneous) by replacing the last {\sf h} in the definition of {\sf J} with {\sf\verb#h^2#} (using \cite{BrSk}), see \cite{DiSt1}, etc.
\msn
{\bf 5.3.~Example.} Let
$$f=(y^2z-x^3)^3+y^9\q\h{with}\q h=(y^2-x^3)^3+y^9,$$
where $n=3$, $d=9$. Here $\mu_Z=58$ (using the above Remark) and $\chi(U)=8\cdot 7+1-58=-1$, see \cite[2.9.5]{wh}. Note that $Z$ is topologically the same as in Example~(5.2).
\sk
Using a computer program based on the algorithm in this paper, we can get
$$\scalebox{0.9}{$\begin{array}{rcccccccccccccccccccccccccccc}
k:&\bk 3 &\bk 4 &\bk 5 &\bk 6 &\bk 7 &\bk 8 &\bk 9 &\bk 10 &\bk 11 &\bk 12 &\bk 13 &\bk 14 &\bk 15 &\bk 16 &\bk 17 &\bk 18 &\bk 19 &\bk 20 &\bk 21 &\bk 22 &\bk 23 &\bk 24 &\bk 25 &\bk 26\\
\gamma_k:&\bk 1 &\bk 3 &\bk 6 &\bk 10 &\bk 15 &\bk 21 &\bk 28 &\bk 36 &\bk 42 &\bk 46 &\bk 48 &\bk 48 &\bk 46 &\bk 42 &\bk 36 &\bk 28 &\bk 21 &\bk 15 &\bk 10 &\bk 6 &\bk 3 &\bk 1 &\\
\mu_k:&\bk 1 &\bk 3 &\bk 6 &\bk 10 &\bk 15 &\bk 21 &\bk 28 &\bk 36 &\bk 42 &\bk 46 &\bk 49 &\bk 51 &\bk 52 &\bk 52 &\bk 52 &\bk 52 &\bk 52 &\bk 52 &\bk 52 &\bk 52 &\bk 52 &\bk 52 &\bk 52 &\bk 52\\
\mu^{\scriptscriptstyle(2)}_k:&\bk 1 &\bk 3 &\bk 4 &\bk 5 &\bk 7 &\bk 7 &\bk 6 &\bk 7 &\bk 7 &\bk 6 &\bk 6 &\bk 6 &\bk 6 &\bk 6 &\bk 6 &\bk 6 &\bk 6 &\bk 6 &\bk 6 &\bk 6 &\bk 6 &\bk 6 &\bk 6 &\bk 6\\
\mu^{\scriptscriptstyle(3)}_k:& &\bk 1 &\bk 1 & &\bk 1 &\bk 1 & &\bk 1 &\bk 1 & & & & & & & & & & & & & &\\
\nu_k:& & & & & & & & & & &\bk 1 &\bk 3 &\bk 6 &\bk 10 &\bk 16 &\bk 24 &\bk 31 &\bk 37 &\bk 42 &\bk 46 &\bk 49 &\bk 51 &\bk 52 &\bk 52\\
\nu^{\scriptscriptstyle(2)}_k:& & & & & & & & & & &\bk 1 &\bk 1 &\bk 1 &\bk 2 &\bk 2 &\bk 2 &\bk 2 &\bk 2 &\bk 2 &\bk 3 &\bk 4 &\bk 5 &\bk 6 &\bk 6\\
\nu^{\scriptscriptstyle(3)}_k:& & & & & & & & & & &\bk 1 &\bk 1 &\bk 1 &\bk 2 &\bk 2 &\bk 2 &\bk 2 &\bk 2 &\bk 1 &\bk 1 &\bk 1 & &\end{array}$}$$
The spectral sequence then degenerates at $E_3$ for $k\les 3d-1$ by Theorem~2, and a variant of condition~(18) with $\R_Z$ replaced by $\{d\}$ does not hold in this case.
\sk
On the other hand, we get by using the program ``bfunction" in RISA/ASIR
$$\tfrac{3}{9}\in\tfrac{1}{9}\,{\rm CS}(f)=[\al_Z,1-\alt_Z]\cap(\R_Z+\Z_{<0})\setminus\R_Z\q\h{and}\q\R_f^0=\emptyset,$$
where $\al_Z=\alt_Z=5/18$. One can show that $\tfrac{3}{9}$ does not belong to the Steenbrink spectral numbers of $h$ (see \cite{St}) by using an embedded resolution of $h^{-1}(0)\subset\C^2$. This does not imply that $\tfrac{3}{9}\notin\R_Z$, although the converse holds (that is, it is a necessary condition). We cannot show that $\tfrac{3}{9}\notin\R_f^0$ by using Theorem~1.
\msn
{\bf 5.4.~Example.} Let
$$f=x^2(x^2+y^2+z^2)^2,$$
where $n=3$, $d=6$. This is {\it non-reduced}, and $\chi(U)=1$. Here $Z_{\rm red}$ is a union of two rational curves intersecting transversally at two points and having degrees 1 and 2. Running a computer program based on the algorithm in this paper, we get
$$\scalebox{0.9}{$\begin{array}{rccccccccccccccccccc}
k:&3 &4 &5 &6 &7 &8 &9 &10 &11 &12 &13 &14 &15 &16 &17 &\cdots\\
\gamma_k:&1 &3 &6 &10 &15 &18 &19 &18 &15 &10 &6 &3 &1\\
\mu_k:&1 &3 &6 &10 &15 &18 &20 &23 &26 &29 &32 &35 &38 &41 &44 &\cdots\\
\mu^{\scriptscriptstyle(2)}_k:&1 &1 &1 &1 &1 &1\\
\nu_k:& & & & & & &1 &5 &11 &19 &26 &32 &38 &44 &50 &\cdots\\
\nu^{\scriptscriptstyle(2)}_k:& & & & & & &1 & & &1\\
\rho_k:& & & & & & & & & & & & &1 &3 &6 &\cdots\\
\rho^{\scriptscriptstyle(2)}_k:& & & & & & & & & & & & &1\end{array}$}$$
Here $\rho_k=\sum_{i=0}^2p^{(2)}_{k-15-i}$ in the notation of (3.5.1). (This is compatible with Corollary~(3.8).)
The calculation implies the $E_2$-degeneration for $k\les 3d-1$ with $p=3$ by Theorem~2. We then get by Corollary~1
$$6\,\R_f^0=\{4,5,7,8\},\q\h{since}\q\R_Z=\bl\{\tfrac{1}{2},1\br\}.$$
This is compatible with a calculation by RISA/ASIR saying that
$$b_f(s)=(s+1)^3\bl(s+\tfrac{1}{2}\br)^3\bl(s+\tfrac{4}{6}\br)\bl(s+\tfrac{5}{6}\br)\bl(s+\tfrac{7}{6}\br)\bl(s+\tfrac{8}{6}\br).$$
\msn
{\bf 5.5.~Example.} Let
$$f=x^4y^2z(x{+}y{+}z)(x{+}y)\q\h{or}\q x^5yz(x{+}y{+}z)(x{+}y),$$
where $n=3$, $d=9$. We have the following for both of the two polynomials:
$$\scalebox{0.88}{$\begin{array}{rcccccccccccccccccccccccccccc}
k:&\bk 3 &\bk 4 &\bk 5 &\bk 6 &\bk 7 &\bk 8 &\bk 9 &\bk 10 &\bk 11 &\bk 12 &\bk 13 &\bk 14 &\bk 15 &\bk 16 &\bk 17 &\bk 18 &\bk 19 &\bk 20 &\bk 21 &\bk 22 &\bk 23 &\bk 24 &\bk 25 &\bk 26\\
\gamma_k:&\bk 1 &\bk 3 &\bk 6 &\bk 10 &\bk 15 &\bk 21 &\bk 28 &\bk 36 &\bk 42 &\bk 46 &\bk 48 &\bk 48 &\bk 46 &\bk 42 &\bk 36 &\bk 28 &\bk 21 &\bk 15 &\bk 10 &\bk 6 &\bk 3 &\bk 1 &\\
\mu_k:&\bk 1 &\bk 3 &\bk 6 &\bk 10 &\bk 15 &\bk 21 &\bk 28 &\bk 36 &\bk 42 &\bk 46 &\bk 50 &\bk 54 &\bk 58 &\bk 62 &\bk 66 &\bk 70 &\bk 74 &\bk 78 &\bk 82 &\bk 86 &\bk 90 &\bk 94 &\bk 98 &\bk 102\\
\mu^{\scriptscriptstyle(2)}_k:&\bk 1 &\bk 1 &\bk 1 &\bk 1 &\bk 1 &\bk 1 &\bk 4 &\bk 1 &\bk 1 &\bk 1\\
\mu^{\scriptscriptstyle(3)}_k:& &\bk 1 &\bk 1 &\bk 1 &\bk 1 &\bk 1 &\bk 4 &\bk 1 &\bk 1 &\bk 1\\
\nu_k:& & & & & & & & & & &\bk 2 &\bk 6 &\bk 12 &\bk 20 &\bk 30 &\bk 42 &\bk 53 &\bk 63 &\bk 72 &\bk 80 &\bk 88 &\bk 96 &\bk 104 &\bk 112\\
\nu^{\scriptscriptstyle(2)}_k:& & & & & & & & & & & & & & & &\bk 4 & & &\bk 1 \\
\nu^{\scriptscriptstyle(3)}_k:& & & & & & & & & & & & & & & &\bk 4\\
\rho_k:& & & & & & & & & & & & & & & & & & & & &\bk 1 &\bk 3 &\bk 6 &\bk 10\\
\rho^{\scriptscriptstyle(2)}_k:\end{array}$}$$
The $\mu_k,\nu_k,\rho_k$ coincide with the ones obtained in Remark~(3.9)(i). The spectral sequence degenerates at $E_3$ for $k\les 3d-1$ by Theorem~2, and we have 
$$9\,\R_f^0=\{3,4,5,7,8,10,11\}$$
since $3/9\in\R_f^0$ by \cite[Theorem 1.4]{Wa2}, although $3\notin{\rm Supp}\,\bl\{\mu_k^{(\infty)}\br\}$. These give examples of non-reduced hyperplane arrangements of 3 variables such that the converse of (2) {\it fails\1} with condition~(2) {\it unsatisfied.}
\msn
{\bf Remark.} The situation is a little bit different in the case
$$f=x^3y^3z(x{+}y{+}z)(x{+}y),$$
where $\R_f^0$ is the same as above, but $\mu^{(3)}_3=1$ so that $3\in\R_f^0$ by combining Theorems~1 and 2. More precisely, we have the following:
$$\scalebox{0.88}{$\begin{array}{rcccccccccccccccccccccccccccc}
k:&\bk 3 &\bk 4 &\bk 5 &\bk 6 &\bk 7 &\bk 8 &\bk 9 &\bk 10 &\bk 11 &\bk 12 &\bk 13 &\bk 14 &\bk 15 &\bk 16 &\bk 17 &\bk 18 &\bk 19 &\bk 20 &\bk 21 &\bk 22 &\bk 23 &\bk 24 &\bk 25 &\bk 26\\
\gamma_k:&\bk 1 &\bk 3 &\bk 6 &\bk 10 &\bk 15 &\bk 21 &\bk 28 &\bk 36 &\bk 42 &\bk 46 &\bk 48 &\bk 48 &\bk 46 &\bk 42 &\bk 36 &\bk 28 &\bk 21 &\bk 15 &\bk 10 &\bk 6 &\bk 3 &\bk 1 &\\
\mu_k:&\bk 1 &\bk 3 &\bk 6 &\bk 10 &\bk 15 &\bk 21 &\bk 28 &\bk 36 &\bk 42 &\bk 46 &\bk 50 &\bk 54 &\bk 58 &\bk 62 &\bk 66 &\bk 70 &\bk 74 &\bk 78 &\bk 82 &\bk 86 &\bk 90 &\bk 94 &\bk 98 &\bk 102\\
\mu^{\scriptscriptstyle(2)}_k:&\bk 1 &\bk 1 &\bk 1 &\bk 2 &\bk 1 &\bk 1 &\bk 4 &\bk 1 &\bk 1 &\bk 1\\
\mu^{\scriptscriptstyle(3)}_k:&\bk 1 &\bk 1 &\bk 1 &\bk 2 &\bk 1 &\bk 1 &\bk 4 &\bk 1 &\bk 1 &\bk 1\\
\nu_k:& & & & & & & & & & &\bk 2 &\bk 6 &\bk 12 &\bk 20 &\bk 30 &\bk 42 &\bk 53 &\bk 63 &\bk 72 &\bk 80 &\bk 88 &\bk 96 &\bk 104 &\bk 112\\
\nu^{\scriptscriptstyle(2)}_k:& & & & & & & & & & & & &\bk 1 & & &\bk 4 & & &\bk 1 \\
\nu^{\scriptscriptstyle(3)}_k:& & & & & & & & & & & & &\bk 1 & & &\bk 4 & & &\bk 1 \\
\rho_k:& & & & & & & & & & & & & & & & & & & & &\bk 1 &\bk 3 &\bk 6 &\bk 10\\
\rho^{\scriptscriptstyle(2)}_k:\end{array}$}$$
Here the $\mu_k,\nu_k,\rho_k$ are the same as above, but the $\mu^{(2)}_k,\nu^{(2)}_k$ are a little bit different, and moreover we have the $E_2$-degeneration (for degree $\les 26$) in this case.
\msn
{\bf 5.6.~Example.} Let
$$f=(x^2+y^2+z^2+w^2)^2+w^4,$$
where $n=d=4$, and $\chi(U)=-2$. Here $Z$ is a union of two smooth quadratic hypersurfaces of $\PP^3$ intersecting each other along a smooth rational curve $C$ with transversal singularity type $A_3$ (that is, locally $\{v^2+u^4=0\}$), and the degree of the curve is 2 so that $\tau_{Z'}=6$ (where $Z'$ is a general hyperplane section of $Z$). Running a computer program based on the algorithm in this paper for $n=4$, we can get
$$\scalebox{0.9}{$\begin{array}{rccccccccccccccccccc}
k:&4 &5 &6 &7 &8 &9 &10 &11 &12 &13 &14 &15 &\cdots\\
\gamma_k:&1 &4 &10 &16 &19 &16 &10 &4 &1\\
\mu_k:&1 &4 &10 &16 &22 &27 &33 &39 &45 &51 &57 &63 &\cdots\\
\mu^{\scriptscriptstyle(2)}_k:&1 &1 &2 &1 &1 \\
\nu_k:& & & & &3 &11 &23 &36 &48 &60 &72 &84 &\cdots\\
\nu^{\scriptscriptstyle(2)}_k:\\
\rho_k:& & & & & & & &1 &4 &9 &15 &21 &\cdots\\
\rho^{\scriptscriptstyle(2)}_k:& & & & & & & &1 &1 &1\end{array}$}$$
Calculating the saturation of the Jacobian ideal $(\dd f)\subset R$ for the maximal ideal $(x,y,z,w)$ by Macaulay2 or Singular as in Remarks~(3.7)(ii) and (4.11)(v) (where $E=(x,y,z,w)$ and $\mu''_k=0$ since $Z$ is analytically locally trivial along the smooth projective curve $C$), and applying Corollary~(3.6), we can get in the notation of (3.5) and (4.2)
$$\scalebox{0.9}{$\aligned&\mu'=v^6{+}v^7{+}v^8,\q\mu''=\mu'''_{\rm def}=0,\q\mu'''=(v^4{+}2\1v^5{+}2\1v^6{+}v^7)\1p^{(2)},\\&\nu'=\nu''=0,\q\rho=(v^{11}{+}2\1v^{12}{+}2\1v^{13}{+}v^{14})\1p^{(2)}.\endaligned$}$$
Here we have to use the above calculation of the $\rho_k$ and apply Corollary~(3.6) to show the vanishing of $\mu'''_{\rm def}$. As is noted in Caution~2 in Remark~(3.7)(ii), however, we {\it cannot\1} decide whether
$$\scalebox{0.9}{$\aligned\bl(\nu'''_{\rm max},\,\nu'''_{\rm def}\br)={}&\bl((4\1v^8{+}4\1v^9{+}4\1v^{10})\1p^{(2)},\,\,v^8{+}v^9{+}v^{10}\br)\\ \h{or}\q&\bl((3\1v^8{+}6\1v^9{+}3\1v^{10})\1p^{(2)},\,\,v^9{+}v^{10}\br),\endaligned$}$$
since this depends on the rank of the middle morphism $\phi$ in the exact sequence (3.4.5).
\sk
By the above calculation together with Theorems~2 and (4.3), the pole order spectral sequence degenerates at $E_2$ for any $k$. By the above description of $Z$, we have
$$\R_Z=\bl\{\tfrac{3}{4},\tfrac{4}{4},\tfrac{5}{4}\br\},\q\h{hence}\q{\rm CS}(f)=\emptyset.$$
So the above computation together with Corollary~1 implies
$$4\,\R_f^0=\{6,7,8\}.$$
\msn
{\bf 5.7.~Example.} Let
$$f=x^3z+x^2y^2+y^3w+y^2w^2+x^2w^2,$$
where $n=d=4$, and $\chi(U)=-1$. This is mentioned in \cite[Example 5.9]{DiSt2} as an example of non-$E_2$-degeneration without any details. The singular locus of $Z$ is a rational curve of degree 1, defined by $x=y=0$, and the generic transversal singularity type is $A_1$ so that $\tau_{Z'}=1$. By a computer program based on the algorithm in this paper for $n=4$, we have
$$\scalebox{0.9}{$\begin{array}{rccccccccccccccccccc}
k:&4 &5 &6 &7 &8 &9 &10 &11 &12 &13 &14 &15 &\cdots\\
\gamma_k:&1 &4 &10 &16 &19 &16 &10 &4 &1\\
\mu_k:&1 &4 &10 &16 &19 &18 &19 &20 &21 &22 &23 &24 &\cdots\\
\mu^{\scriptscriptstyle(2)}_k:&1 &2 &2 &2 &2 &1 &1 &1 &1 &1 &1 &1 &\cdots \\
\mu^{\scriptscriptstyle(3)}_k:&1 &1 &1 &1 &1\\
\nu_k:& & & & & &2 &9 &16 &20 &22 &24 &26 &\cdots\\
\nu^{\scriptscriptstyle(2)}_k:& & & & & & & & & &1 &1 &1 &\cdots\\
\nu^{\scriptscriptstyle(3)}_k:\\
\rho_k:& & & & & & & & & & &1 &2 &\cdots\\
\rho^{\scriptscriptstyle(2)}_k:\end{array}$}$$
Calculating the saturation of the Jacobian ideal $(\dd f)\subset R$ for the maximal ideal $(x,y,z,w)$ and for the ideal $E=(x,y,w)$ by Macaulay2 or Singular as in Remarks~(3.7)(ii) and (4.11)(v), and applying Corollary~(3.6), we get in the notation of (3.5) and (4.2)
$$\scalebox{0.9}{$\aligned&\mu'=2\1v^8,\,\,\,\mu''_k=(2\1v^5{+}5\1\1v^6{+}5\1\1v^7)\1p^{(1)},\,\,\,\mu'''_k=v^4\1p^{(2)},\,\,\,\mu'''_{\rm def}=0,\\&\nu'=0,\q\nu''=(5\1\1v^{10}{+}5\1\1v^{11}{+}2\1\1v^{12})\1p^{(1)},\q\rho_k=v^{14}\1p^{(2)}.\endaligned$}$$
Here the vanishing of $\mu'''_{{\rm def},\,k}$ follows from $\mu'''_k=p^{(2)}_{k-4}$, and the above data for $\rho_k$ is not needed (since the rank is 1 in this case). As is noted in Caution~2 in Remark~(3.7)(ii), however, we {\it cannot\1} decide whether
$$\scalebox{0.9}{$\bl(\nu'''_{\rm max},\,\nu'''_{\rm def}\br)=\bl((v^8{+}v^{10})\1p^{(2)},\,v^8\br)\,\,\,\,\h{or}\,\,\,\,\bl(2\1v^9\1p^{(2)},\,0\br),$}$$
since this depends on the rank of the morphism $\phi$ in (3.4.5). Here we see that the rank is 1 or 2, but not 0 (since we get a contradiction if it is 0 and $\nu'''_{{\rm def},\,8}=2$).
\sk
By the above calculation and Theorem~2, the pole order spectral sequence degenerates at $E_3$ for $k\les 4d-1$. On the other hand, we get by RISA/ASIR
$$4\,\R_Z\cap\Z=\{4,6\},\q\al_Z=5/6,\q\h{hence}\q{\rm CS}(f)=\emptyset.$$
The above calculation together with Corollary~1 then implies
$$4\,\R_f^0=\{5,7,8\}.$$
\msn
{\bf 5.8.~Example.} Let
$$f=xyzw(x+y+z)(y-z+w),$$
where $n=4$, $d=6$, and $\chi(U)=-2$. This is treated in \cite[Example 5.7]{DiSt2}. The non-free locus of $Z$ consists of two points $z_i$ ($i=1,2$) such that $(Z,z_i)$ is locally isomorphic to a generic central hyperplane arrangement of $(\C^3,0)$ with degree 4. The generic hyperplane section of $Z$ has only $A_1$-singularities, and there are 15 singular points so that $\tau_{Z'}=15$. By a computer program based on the algorithm in this paper for $n=4$, we can get
$$\scalebox{0.9}{$\begin{array}{rccccccccccccccccccccccccc}
k:&\bk 4 &\bk 5 &\bk 6 &\bk 7 &\bk 8 &\bk 9 &\bk 10 &\bk 11 &\bk 12 &\bk 13 &\bk 14 &\bk 15 &\bk 16 &\bk 17 &\bk 18 &\bk 19 &\bk 20 &\bk 21 &\bk 22 &\bk 23 &\bk\cdots\\
\gamma_k:&\bk 1 &\bk 4 &\bk 10 &\bk 20 &\bk 35 &\bk 52 &\bk 68 &\bk 80 &\bk 85 &\bk 80 &\bk 68 &\bk 52 &\bk 35 &\bk 20 &\bk 10 &\bk 4 &\bk 1 & & &\\
\mu_k:&\bk 1 &\bk 4 &\bk 10 &\bk 20 &\bk 35 &\bk 52 &\bk 68 &\bk 82 &\bk 97 &\bk 112 &\bk 127 &\bk 142 &\bk 157 &\bk 172 &\bk\cdots\\
\mu^{\scriptscriptstyle(2)}_k:&\bk 1 &\bk 2 &\bk 8 &\bk 2 &\bk 2 &\bk 2 &\bk 1 & & & & & & &\\
\nu_k:& & & & & & & &\bk 2 &\bk 12 &\bk 32 &\bk 59 &\bk 90 &\bk 122 &\bk 152 &\bk\cdots\\
\nu^{\scriptscriptstyle(2)}_k:& & & & & & & & &\bk 10 & & & & &\\
\rho_k:& & & & & & & & & & & & & & &5 &\bk 14 &\bk 26 &\bk 40 &\bk 55 &\bk 70 &\bk\cdots\\
\rho^{\scriptscriptstyle(2)}_k:& & & & & & & & & & & & & & &5 & & & & &\\
\end{array}$}$$
Calculating the saturation of the Jacobian ideal $(\dd f)\subset R$ for the maximal ideal $(x,y,z,w)$ and for the ideal $E=(y,z,xw)$ by Macaulay2 or Singular as in Remarks~(3.7)(ii) and (4.11)(v), and applying Corollary~(3.6), we can get in the notation of (3.5) and (4.2)
$$\scalebox{0.9}{$\aligned&\mu'=v^{10},\q\mu''=2\1v^9\1p^{(1)},\q\mu'''_{\rm def}=0,\q\mu'''=\msum_{i=0}^4\,(i+1)\1v^{i+4}\1p^{(2)},\\&\nu'=0,\q\nu''=2\1v^{16}\1p^{(1)},\q\rho=\msum_{i=0}^4\,(5-i)\1v^{i+18}\,p^{(2)},\\&\nu'''_{\rm max}=\bl(2\1v^{11}{+}8\1v^{12}{+}10\1v^{13}{+}8\1v^{14}{+}2\1v^{15}\br)\1p^{(2)},\q\nu'''_{\rm def}=v^{14}.\endaligned$}$$
(The non-vanishing of $\mu''_k$ means that $Z$ is not a free divisor.) Here we have to use the above calculation of the $\rho_k$ and apply Corollary~(3.6) to show the vanishing of the $\mu'''_{{\rm def},\,k}$. There is no problem related to the morphism $\phi$ in (3.4.5), since $\mu'_k=0$ except for $k=10\,(\ne 2d)$ and (3.7.3) holds. (Related to Remark~(3.9)(ii), we can construct a filtration on $(M''')^{\sim}_{\PP^3}$ such that the graded pieces of its direct image by $\overline{\pi}_2$ give the direct factors although (SC) is {\it not\1} satisfied.)
\sk
By the above calculation together with Theorem~(4.3), the pole order spectral sequence degenerates at $E_2$ for {\it any\1} $k$. (Here we have to use the above calculation of the $\mu'_k$, etc.\ in order to show that $\chi_{f,k}=0$ for {\it any} $k\ges 2d-1$.) This is compatible with a computation of the $\mu_k^{(2)}$ in \cite{DiSt2}, and the latter can be justified by this calculation.
\sk
By the above description of $Z$ together with \cite{Wa2} (or \cite{bha}) we have
$$\R_Z=\bl\{\tfrac{3}{4},\tfrac{4}{4},\tfrac{5}{4},\tfrac{6}{4}\br\},\q\h{hence}\q{\rm CS}(f)=\emptyset.$$
This coincides with a calculation by RISA/ASIR. So the above computation together with Corollary~1 gives
$$6\,\R_f^0=\{4,5,7,8,10\}.$$
\sk
This is one of rather nontrivial examples of central hyperplane arrangements in $\C^4$ with degree 6. One may also calculate the case of $xyzw(x+y+z)(z+w)$ with $\chi(U)=-1$.
\msn
{\bf 5.9.~Example.} Let
$$f=x^6+x^4yz+y^3w^3+y^6,$$
where $n=4$, $d=6$, and $\chi(U)=-1$. General hyperplane sections of $Z$ have a singularity defined locally analytically by $u^6+cu^4v+v^3=0$ ($c\in\C$) using $\C^*$-action and \cite{KaSc}, \cite{Tj}, \cite{Va2}. This is one of examples with condition~(AT) {\it unsatisfied}, see also Remark~(3.7)(iii). By a computer program based on the algorithm in this paper for $n=4$, we can get
$$\scalebox{0.9}{$\begin{array}{rccccccccccccccccccccccccc}
k:&\bk 4 &\bk 5 &\bk 6 &\bk 7 &\bk 8 &\bk 9 &\bk 10 &\bk 11 &\bk 12 &\bk 13 &\bk 14 &\bk 15 &\bk 16 &\bk 17 &\bk 18 &\bk 19 &\bk 20 &\bk 21 &\bk 22\\
\gamma_k:&\bk 1 &\bk 4 &\bk 10 &\bk 20 &\bk 35 &\bk 52 &\bk 68 &\bk 80 &\bk 85 &\bk 80 &\bk 68 &\bk 52 &\bk 35 &\bk 20 &\bk 10 &\bk 4 &\bk 1\\
\mu_k:&\bk 1 &\bk 4 &\bk 10 &\bk 20 &\bk 35 &\bk 52 &\bk 68 &\bk 81 &\bk 90 &\bk 98 &\bk 106 &\bk 115 &\bk 124 &\bk 133 &\bk 142 &\bk 151 &\bk 160 &\bk 169 &\bk\cdots\\
\mu^{\scriptscriptstyle(2)}_k:&\bk 1 &\bk 3 &\bk 5 &\bk 5 &\bk 4 &\bk 4 &\bk 4 &\bk 4 &\bk 3 &\bk 2 &\bk 2 &\bk 2 &\bk 2 &\bk 2 &\bk 2 &\bk 2 &\bk 2 &\bk 2 &\bk\cdots\\
\mu^{\scriptscriptstyle(3)}_k:& & &\bk 1 &\bk 2 &\bk 2 &\bk 2 &\bk 2 &\bk 2 &\bk 1\\
\nu_k:& & & & & & & &\bk 1 &\bk 5 &\bk 18 &\bk 38 &\bk 63 &\bk 89 &\bk 113 &\bk 132 &\bk 149 &\bk 165 &\bk 183 &\bk\cdots\\
\nu^{\scriptscriptstyle(2)}_k:& & & & & & & & & &\bk 1 &\bk 1 &\bk 1 &\bk 2 &\bk 4 &\bk 4 &\bk 3 &\bk 2 &\bk 2 &\bk\cdots\\
\nu^{\scriptscriptstyle(3)}_k:& & & & & & & & & &\bk 1 &\bk 1 &\bk 1 &\bk 1 &\bk 1 & & & &\\
\rho_k:& & & & & & & & & & & & & & & &\bk 2 &\bk 6 &\bk 14 &\bk 23\\
\rho^{\scriptscriptstyle(2)}_k:\end{array}$}$$
This computation takes very long, perhaps more than 10 hours, depending on the computer, if we verify $\mu^{(3)}_k=0$ until $k=21$. Calculating the saturation of the Jacobian ideal $(\dd f)\subset R$ for the maximal ideal $(x,y,z,w)$ and for the ideal $E=(x,y,w)$ by Macaulay2 or Singular as in Remarks~(3.7)(ii) and a variant of Remark~(4.11)(v), and applying Corollary~(3.6), we can get in the notation of (3.5) and (4.2)
$$\scalebox{0.9}{$\aligned\mu'&=v^9{+}2\1v^{10}{+}3\1v^{11}{+}2\1v^{12}{+}v^{13},\\
\mu''&=(v^7{+}4\1v^8{+}6\1v^9{+}6\1v^{10}{+}3\1v^{11}{+}v^{12})\1p^{(1)},\\
\mu'''_{\rm max}&=(v^4{+}4\1v^5{+}2\1v^6{+}2\1v^7)\1p^{(2)},\\
\mu'''_{\rm def}&=2\1v^5{+}3\1v^6{+}3\1v^7{+}v^8,\\
\nu'&=v^{16}{+}3\1v^{17}{+}3\1v^{18}{+}2\1v^{19},\\
\nu''&=(v^{13}{+}3\1v^{14}{+}6\1v^{15}{+}6\1v^{16}{+}4\1v^{17}{+}v^{18})\1p^{(1)},\\
\rho&=(2\1v^{19}{+}2\1v^{20}{+}4\1v^{21}{+}v^{22})\1p^{(2)},\endaligned$}$$
However, we {\it cannot\1} determine the integers $a\in[0,2]$, $b\in[0,1]$ satisfying
$$\scalebox{0.9}{$\aligned\nu'''_{\rm max}&=(2\1v^{11}{+}3\1v^{12}{+}8\1v^{13}{+}3\1v^{14}{+}2\1v^{15})\1p^{(2)}\\&\q-a\1(1-v)^2v^{12}\1p^{(2)}-b\1(1-v)^2(v^{11}+v^{13})\1p^{(2)},\\ \nu'''_{\rm def}&=v^{11}{+}2\1v^{12}{+}3\1v^{13}{+}2\1v^{14}{+}v^{15}\\&\q-a\1v^{12}-b\1(v^{11}+v^{13}).\endaligned$}$$
\sk
By the above calculation and Theorem~2, the pole order spectral sequence degenerates at $E_3$ for $k\les 4d-1$. Although it is not easy to determine $\R_Z$ for this example, we can show the following by using the above calculation together with Theorem~1:
$$6\1\R_f^0\cup(6\1\R_Z\cap\Z)=\{9,\dots,12\}\cup(6\1\R_Z\cap\Z),$$
since $6\1\R_Z\supset\{3,\dots,8\}$ by computing $b_h(s)$ for $h=f|_{w=1}$ using RISA/ASIR (and $\R_Z\subset(0,3)$ by \cite{mic}).
\msn
{\bf Remark.} There are small computer programs for $n\,{=}\,3,4$, based on the algorithm in this paper, and running on MacOS (probably 10.9.5 or later) and on Unix (SPARC/Solaris 10).
\bs\bs
\vbox{
\centerline{\bf Appendix. Double symmetry of modified spectra}
\bsn
In this Appendix we show a {\it double symmetry\1} of the modified pole-order spectrum for the $E_2$-term in the case of strongly free divisors on $\PP^3$.}
\msn
{\bf A.1.} Let $Z$ be a strongly free reduced divisor on $\PP^3$ (see (4.1)) with $f$ a defining polynomial. Put $d=\deg f$. Assume $Z$ {\it essential,} that is, $f$ is not a polynomial of 3 variables. We have the {\it modified pole-order spectrum for the $E_2$-term} as is defined in (4.2):
$$\chi_f:=\mu-v^{-d}\nu+v^{-2d}\rho\,\in\,\Z[v].$$
It has the following symmetry as is shown in Theorem~(4.5):
$$v^{2d}\chi_f(v^{-1})=\chi_f(v).
\leqno{\rm(A.1.1)}$$
It turns out, however, that there is a {\it hidden second symmetry\1} as follows.
\sk
Set
$$\chih_f:=\chi_f-v^d,\q{}^D\chih_f:=(1-v)\chih_f,$$
\vskip-3mm
$${}^D\chih_f^{\,\les d}:=\msum_{k\les d}\,{}^D\chih_{f,k}\1v^k,\q{}^D\chih_f^{\,>d}:={}^D\chih_f-{}^D\chih_f^{\,\les d},$$
with
$${}^D\chih_f=\msum_{k\in\Z}\,{}^D\chih_{f,k}\1v^k,\q{}^D\chih_{f,k}=\chih_{f,k}-\chih_{f,k-1},\q\chih_{f,k}=\chi_{f,k}-\delta_{k,d}.$$
Note that $^D$ corresponds to ${\rm Diff}$ in (3.5).
\msn
{\bf A.2.~Theorem.} {\it We have the following anti-symmetry and symmetry}\1:
$${}^D\chih_f^{\,\les d}(v^{-1})v^{2d+1}=-{}^D\chih_f^{\,>d}(v),
\leqno{\rm(A.2.1)}$$
$${}^D\chih_f^{\,\les d}(v^{-1})v^{d+1}={}^D\chih_f^{\,\les d}(v).
\leqno{\rm(A.2.2)}$$
\msn
{\bf Remark.} The anti-symmetry (A.2.1) is equivalent to the symmetry (A.1.1), where Diff transforms a symmetry into an anti-symmetry. The symmetry of ${}^D\chih_f^{\,\les d}$ in (A.2.2) does not seem to be noticed very well. 
\msn
{\bf A.3.~Proof of Theorem~A.2.} As in the proof of Theorem~(4.5), set in $\Z[v]$
$$\chit_f:=(1-v)^2\chi_f,\q\mut:=(1-v)^2\mu,\q\nut:=(1-v)^2\nu,\q\rhot:=(1-v)^2\rho.$$
Put
$$\aligned&\q\q\q\gamma\1':=\dfrac{(v^d-v)^2}{(v-1)^2}=\bl(v+\cdots+v^{d-1}\br)^2\\={}&v^2+2v^3+\cdots+(d-1)v^d+\cdots+2v^{2d-3}+v^{2d-2}.\endaligned$$
Since
$$\mu-\nu+\rho=\gamma\1'{}^2,$$
we have
$$\aligned\mut-\nut+\rhot&=(v^d-v)^2\gamma\1'\\&=v^2\gamma\1'-2v^{d+1}\gamma\1'+v^{2d}\gamma\1'.\endaligned$$
Set
$$\at:=v^2\gamma\1'-\mut,\q \bt:=2v^{d+1}\gamma\1'-\nut,\q \ct:=v^{2d}\gamma\1'-\rhot\q\h{in}\q\Z[v],$$
so that
$$\at-\bt+\ct=0.
\leqno{\rm(A.3.1)}$$
\sk
From Corollary~(3.6) and (4.1.2), we can deduce that
$$v^{4d+2}\mut(v^{-1})=\rhot(v),\q v^{4d+2}\nut(v^{-1})=\nut(v).$$
Since
$$v^{2d}\gamma\1'(v^{-1})=\gamma'(v),$$
we then get
$$v^{4d+2}\at(v^{-1})=\ct(v),\q v^{4d+2}\bt(v^{-1})=\bt(v).
\leqno{\rm(A.3.2)}$$
\sk
By Lemma~(A.4) below, we have
$${\rm Supp}\bl\{\at_k\br\}\subset[d+3,2d],\q {\rm Supp}\bl\{\ct_k\br\}\subset[2d+2,3d-1],
\leqno{\rm(A.3.3)}$$
with
$$\at_{2d}=\ct_{2d+2}=1.
\leqno{\rm(A.3.4)}$$
Here we use the assumption that $Z$ is essential.
\sk
We then see that
$$\aligned\chit_f&=\mut-v^{-d}\nut+v^{-2d}\rhot\\&=(1-v)^2\gamma\1'-\at+v^{-d}\1\bt-v^{-2d}\ct\\&=v^{2d}-2v^{d+1}+v^2+(v^{-d}-1)(\at-v^{-d}\ct),\endaligned
\leqno{\rm(A.3.5)}$$
with
$$\chit_{f,2}=0,\q\chit_{f,d}=1,\q\chit_{f,d+1}=-2,\q\chit_{f,d+2}=1,\q\chit_{f,2d}=0.
\leqno{\rm(A.3.6)}$$
Set
$$\varepsilon=v^{-d}(\at-v^{-d}\ct)+v^2-v^d.$$
From the above assertions we can deduce that
$$\chit_f-(1-v)^2v^d=\varepsilon-v^d\varepsilon,
\leqno{\rm(A.3.7)}$$
$$v^{d+2}\varepsilon(v^{-1})=-\varepsilon,\q{\rm Supp}\bl\{\varepsilon_k\br\}\subset[3,d-1].
\leqno{\rm(A.3.8)}$$
These imply the assertions (A.2.1--2). This finishes the proof of Theorem~(A.2).
\msn
{\bf A.4.~Lemma.} {\it In the above notation and assumption, the assertions {\rm (A.3.3--4)} hold.}
\msn
{\it Proof.} For (A.3.3) it is enough to show the following inclusion in view of (A.3.2):
$${\rm Supp}\bl\{\at_k\br\}\subset[d{+}3,+\infty).
\leqno{\rm (A.4.1)}$$
Indeed, the inclusion ${\rm Supp}\bl\{\ct_k\br\}\subset[2d{+}2,+\infty)$ holds, since we have by definition
$$\rho_k=\rhot_k=0\,\,\,(\forall\,k<2d+2).
\leqno{\rm(A.4.2)}$$
\sk
By the definition of $\at_k$, the assertion (A.4.1) is equivalent to
$$\mut_k=k-3\,\,\,\,(\forall\,k\in[4,d+2]).
\leqno{\rm(A.4.3)}$$
Set
$$M_{\les j}:=\mopl_{k\les j}\,M_k,\,\,\,\h{and similarly for}\,\,\,\,R_{\les j}.$$
Then
$$M_{\les d+2}(4)=R_{\les d-2},
\leqno{\rm(A.4.4)}$$
since $\deg f_i=d-1$ with $f_i:=\dd f/\dd x_i$. Here $x_1,\dots,x_4$ are the coordinates of $\C^4$, and $(4)$ on the left-hand side denotes the shift of grading by $4$. We may assume that $x_1,x_2$ are the pull-back of the coordinates of $\C^2$ by the projection $\pi_2:\C^4\to\C^2$. Then (A.4.3) follows from the graded isomorphism
$$R=\C[x_1,x_2]\otimes_{\C}\C[x_3,x_4].$$
\sk
As for (A.3.4), the assertion is equivalent to the following by the definition of $\ct_k$ (and using (A.3.2), (A.4.2)):
$$\rho_{2d+2}=0
\leqno{\rm(A.4.5)}$$
So (A.3.4) is reduced to that
$$\om\in\Om^2_2\,\,\,\h{vanishes if}\,\,\,\df\sw\om=0,
\leqno{\rm(A.4.6)}$$
and the last assertion follows from the essentiality of $Z$. Indeed, any nonzero $\om\in\Om^2_2$ is a $\C$-linear combination of $\ddd x_i\sw\ddd x_j$ ($i<j$), and the vanishing of $\df\sw\om$ implies a nontrivial $\C$-linear relation between the partial derivatives $f_i$. However, there is no such relation by the essentiality assumption. (If there is such a relation, then we have the vanishing of $f_i$ for some $i$ after a $\C$-linear coordinate change of $\C^4$, but this contradicts the essentiality of $Z$.) This finishes the proof of Lemma~(A.4).
\msn
{\bf A.5.~Question.} $\q\at_k,\bt_k,\ct_k\ges 0\,\,(k\in\Z),\q\varepsilon_k\ges 0\,\,(k\les d/2)\,\,?$
\msn
This may be related to condition~(18) in the introduction.
\bsn
{\bf A.6.~Remark.} We can verify (A.2.1--2) by using Macaulay2, for instance, as follows:
\ms\vbox{\small\sf\verb#R=QQ[x,y,z,w]; f=x*y*z*w*(x+y+z)*(z+w)*(x+z)*(x+y);#
\sk
\verb#d=first degree f; A=QQ[v]/(v^(4*d)); B=QQ[v]/(v^(d+1)); vA=sub(v,A);#
\sk
\verb#muA=vA^4*sub(hilbertSeries(R/ideal(jacobian ideal(f)),Order=>4*d),vars A);#
\sk
\verb#F=frac(QQ[v]); seqA=sub((1-v^(4*d))/(1-v),A); DDmu=sub(muA*(1-vA)^2,F);#
\sk
\verb#DDrho=sub(DDmu,v=>1/v)*v^(4*d+2); rhoA=sub(DDrho,A)*seqA^2;#
\sk
\verb#nuA=muA+rhoA-sub(((v^d-v)/(v-1))^4,A); DDnu=sub(nuA*(1-vA)^2,F);#
\sk
\verb#rhoSSA=sub(DDrho/v^(2*d),A)*seqA^2; nuSA=sub(DDnu/v^d,A)*seqA^2;#
\sk
\verb#Ch=sub(muA-nuSA+rhoSSA,F); DChh=(Ch-v^d)*(1-v)#
\sk
\verb#DChh1=sub(sub(DChh,B),F); sub(DChh1,v=>1/v)*v^(d+1)-DChh1#
\sk
\verb#DChh2=DChh-DChh1; sub(DChh1,v=>1/v)*v^(2*d+1)+DChh2#}
\msn
If the final two outputs are both zero, then these confirm (A.2.2) and (A.2.1) respectively. We can verify whether a reduced divisor is strongly free by using {\small\sf size(minbase(syz(jacob(f))))} in Singular, for instance. The divisor is strongly free if the output is 3 in the 4 variable case.

\end{document}